\documentclass[11pt]{article}
\usepackage{graphicx}
\usepackage{amsmath,amssymb,amsthm}
\usepackage{mathtools}
\usepackage[table]{xcolor}
\usepackage{caption}
\usepackage{subcaption}
\usepackage{comment}
\usepackage{dsfont}
\usepackage{nicematrix}
\usepackage{enumitem}
\usepackage{authblk}
\usepackage{kotex}
\usepackage{mathrsfs}
\usepackage{overpic}
\usepackage{tikz}
\usepackage{makecell}
\usepackage{hyperref}
\hypersetup{
	colorlinks=true,
	breaklinks=true,
	linkcolor=black,
	urlcolor=blue,
	anchorcolor=blue,
	citecolor=blue
}
\usepackage[nameinlink,capitalise,noabbrev]{cleveref}

\newcommand{\bx}{\mbox{\boldmath $x$}}
\newcommand{\bv}{\mbox{\boldmath $v$}}
\newcommand{\bw}{\mbox{\boldmath $w$}}

\newcommand{\sX}{\mathrm X}
\newcommand{\sx}{\mathrm x}
\newcommand{\sY}{\mathrm Y}

\newtheorem{theorem}{Theorem}[section]
\newtheorem{definition}[theorem]{Definition}
\newtheorem{lemma}[theorem]{Lemma}  
\newtheorem{proposition}[theorem]{Proposition}
\newtheorem{corollary}[theorem]{Corollary}
\newtheorem{remark}[theorem]{Remark}

\title{Adaptive Cucker-Smale Networks:\\ Limiting Laplacian Time-Varying Dynamics}
\author[1,2]{Christian Kuehn}
\author[3]{Jaeyoung Yoon}
\affil[1]{Technical University of Munich, School of Computation Information and Technology, Department of Mathematics, email: ckuehn@ma.tum.de}
\affil[2]{Complexity Science Hub Vienna, Vienna, Austria}
\affil[3]{Technical University of Munich, School of Computation Information and Technology, Department of Mathematics,email: jaeyoung.yoon@tum.de, wodud1516@gmail.com}

\date{\today}

\begin{document}

\maketitle

\begin{abstract}
    Differences in opinion can be seen as distances between individuals, and such differences do not always vanish over time. In this paper, we propose a modeling framework that captures the formation of opinion clusters, based on extensions of the Cucker–Smale and Hegselmann–Krause models to a combined adaptive (or co-evolutionary) network. Reducing our model to a singular limit of fast adaptation, we mathematically analyze the asymptotic behavior of the resulting Laplacian dynamics over various classes of temporal graphs and use these results to explain the behavior of the original proposed adaptive model for fast adaptation. In particular, our approach provides a general methodology for analyzing linear consensus models over time-varying networks that naturally arise as singular limits in many adaptive network models.
\end{abstract}

\section{Introduction} 

Adaptive dynamical networks are models in which the network structure evolves in response to the dynamic states of its nodes, resulting in the co-evolution of both network topology and node behavior. Unlike temporal networks, where the structure is predefined and evolves according to fixed rules, adaptive dynamical networks enable real-time structural changes, offering greater flexibility in modeling complex, real-world systems.

A key advantage of adaptive dynamical networks lies in their ability to model systems where interactions between nodes are not static over time \cite{BGKKY2023,BSY2019,BVSY2021,FSM2021,KYSN2017,Kuehn2015}. In neuroscience, for example, synaptic plasticity describes how neural connections evolve with repeated stimulation or experience, making it an ideal candidate for adaptive network models. Similarly, in social systems, opinion dynamics can be captured by networks that adjust their structure as individuals interact, form, or dissolve groups, providing a more accurate depiction of opinion convergence, polarization, and multi-clustering phenomena.

One example of adaptive dynamical networks is the adaptive Kuramoto (KM) model, which has been widely applied to explain collective behaviors of coupled oscillators such as neural networks \cite{ABMH2022,MLHBT2007,PYT2013,TBTSY2023}. In the context of synaptic plasticity, this model is adapted to reflect how neural synchronization evolves over time as synaptic strengths adjust. This mirrors the biological process of synaptic plasticity, where the strength of connections between neurons changes based on the temporal relationship of their firing patterns. The dynamic evolution of network structure is crucial for understanding how brain networks reorganize in response to learning and experience.

In this paper, we introduce the application of adaptive dynamical networks to opinion dynamics. The Cucker-Smale (CS) model, a foundational framework in opinion dynamics, models agents who update their opinions based on interactions with neighbors in Euclidean space. It differs from the KM regime in that it applies to nodes with unbounded phase spaces, usually in dimensions greater than one, whereas the KM model typically has node phase spaces in finite domains such as the circle $\mathbb S^1$. The classical CS model predicts asymptotic consensus among agents when the interaction kernel is non-negative. However, in real-world scenarios, full consensus is uncommon, as individuals with similar opinions tend to reinforce each other, while those with opposing views often experience friction or exclusion.

To address this limitation, we propose an adaptive CS model in which the network structure evolves as individuals' opinions change, based on the neighbor selection and normalization principles of the Hegselmann–Krause model \cite{AC2018,BHT2010,CF2012,HK2002,PR2021}, albeit in a modified form. This extension reflects the dynamic nature of real-world opinion formation, where the network of opinions adapts over time in response to interactions among individuals. By incorporating adaptive network structures, our model provides a more realistic and subtle representation of opinion dynamics, capturing both the convergence of opinions and the emergence of subgroups within a population.

From a mathematical perspective, this paper focuses on the asymptotic behavior of the adaptive CS model in the singular limit ($\varepsilon\to\infty$) of fast adaptation. This limit then leads to a network dynamical system, which we study using Laplacian dynamics on time-dependent graphs. Such Laplacian dynamics have been extensively studied in the literature over the past decades \cite{CI2012,DaCunha2005,Gallier2016,JKL2024,KSLLDA2010,OSM2004,PSW2003,VL2020,ZC2014}. The novelty of this work lies in establishing lower bounds on the eigenvalues of Laplacian matrices associated with time-continuously evolving complete graphs, enabling a spectral approach to the convergence analysis. Moreover, we establish sufficient connectivity conditions for convergence, and quantify the extent of admissible asymmetry in the interaction structure that preserve the convergence toward a consensus.

\subsection{Derivation of adaptive Cucker-Smale model}\label{model_deri}

Our modeling approach originates from the question of whether the possibilities of dynamic networks, demonstrated through the KM model over the past several decades, can also be applied to Newtonian mechanical motion. To address this question, we can investigate it using the CS model, which describes the physical motion and synchronization of a group. The classical CS model is a semi-linear consensus model in the form of a Laplacian, combined with a distance-based kernel: for position and velocity vectors $\bx_i,\bv_i\in\mathbb R^n$ of the $i$-th particle,
\begin{align*}
    \dot\bx_i=\bv_i,\quad\dot\bv_i=\frac{\kappa}{N}\sum_{j=1}^N\phi(\|\bx_i-\bx_j\|)(\bx_j-\bx_i),\quad\forall~i=\{1,\cdots,N\}=:[N].
\end{align*}
In particular, since the kernel function $\phi$ is typically assumed to be positive, it induces strong velocity alignment among nearby particles according to this distance-based function, and this synchronization (in velocity) gradually spreads to all particles. This represents the classical dynamical scenario of the CS model. Our idea is to interpret this kernel as a time-varying weighted graph. By extending it to allow for negative weights, we propose that a different scenario could arise. Specifically, using $\kappa_{ij}$ to represent the weights (or individual coupling strengths), we consider the following adaptive CS model:
\begin{align*}
    \begin{dcases}
        \dot\bx_i=\bv_i,\quad\dot\bv_i=\frac{\kappa}{N}\sum_{j=1}^N\kappa_{ij}(\bv_j-\bv_i),\\
        \dot\kappa_{ij}=\varepsilon f_{ij}(\bx,\bv,\boldsymbol{\kappa}),
    \end{dcases}
\end{align*}
where $\bx$, $\bv$ and $\boldsymbol{\kappa}$ stands for collections of position, velocity and individual coupling coefficients of $N$ particles. Nonnegative constants $\kappa$ and $\varepsilon$ are a global coupling strength and an adaptation time constant, respectively.

The problem for determining the adaptive rule for $\kappa_{ij}$ is somewhat challenging. The KM model, which has been extensively studied in the context of dynamics networks, is fundamentally based on bounded interaction term involving trigonometric functions. As a result, the model typically works well for any sufficiently smooth adaptive rule. In contrast, the CS model, where velocity changes are governed by a linear combination of relative velocities, can lead to excessive increases in kinetic energy in a short time if an inappropriate adaptive rule is applied.

With this concern in mind, we aim to model an adaptive network that incorporates the concept of bounded confidence as a distinguishing feature from the original CS model. Bounded confidence refers to the psychological and sociological phenomenon where individuals are inclined to disregard opinions significantly different from their own and engage only with those holding similar views. This concept is rooted in confirmation bias and social trust, reflecting the difficulty in persuasion or influence as opinion disparity increases. It was mathematically formalized in the Hegselmann-Krause (HK) model in \cite{HK2002} as follows:
\begin{align*}
    x_i(t+1)&=\frac{1}{|N_i(t)|}\sum_{j\in N_i(t)}x_j(t),\quad\mbox{where}\quad
    N_i(t)=\{j\in[N]~:~|x_i(t)-x_j(t)|<\varepsilon\}.
\end{align*}
Inspired by choosing neighbors and normalization in HK model, we consider the adaptive CS model as follows:
\begin{align}\label{main_model}
    \begin{aligned}
        \begin{dcases}
            \dot\bx_i=\bv_i,\quad\dot\bv_i=\frac{\kappa}{\mathcal I_i(t)}\sum_{j=1}^N\psi_{ij}(t)\kappa_{ij}(t)(\bv_j-\bv_i),\\
            \dot\kappa_{ij}=\varepsilon\left(\frac{\bv_i\cdot\bv_j}{|\bv_i||\bv_j|}-\kappa_{ij}\right),\quad\forall~i,j\in[N].
        \end{dcases}
    \end{aligned}
\end{align}
In other words, it is a CS model on a temporal graph defined as
\begin{align*}
    \left(\frac{\kappa}{\mathcal I_i(t)}\psi_{ij}(t)\kappa_{ij}(t)\right)_{ij},
\end{align*}
where normalization factors $\mathcal I_i(t)$ and connection variables $\psi_{ij}(t)$ are defined by
\begin{align}\label{I_Psi}
    \begin{aligned}
        \mathcal I_i(t)&:=\big|\{j\in[N]~:~\|\bx_i(t)-\bx_j(t)\|<d\}\big|\quad\big(\Rightarrow\mathcal I_i\ge1\big)\\
        \psi_{ij}(t)&:=\chi_{[0,d)}(\|\bx_i(t)-\bx_j(t)\|),
    \end{aligned}
\end{align}
for a hyperparmeter radius $d>0$. Here, $\mathcal I_i$ and $\psi_{ij}$ represent the normalization and connection to neighbors within a radius $d$. The adaptive rule ensures that $\kappa_{ij}$ remains constrained within the interval $[-1,1]$ whenever $\kappa_{ij}^0\in[-1,1]$, thereby preventing excessive energy flowing into the system. Furthermore, we can describe the influence magnitude, which becomes stronger/weaker depending on opinion similarity/disparity, which the HK model does not capture explicitly.

The key distinction lies in that the HK model is a single-layer framework focusing solely on the opinions of particles, whereas, our approach employs a dual-layer model where position and velocity vectors are treated as spatial coordinate and opinion, respectively. This dual-layer structure implies that physical distances may also change according to opinion differences \cite{BHOY2022}. One significant mathematical consequence of this type of neighbor interaction and normalization is the discontinuity in the vector field, which undermines the existence and uniqueness of classical solutions. This issue will be discussed in \Cref{subsec_existence}.

In summary, our adaptive modeling framework is a combination of Cucker-Smale, adaptive Kuramoto, and Hegselmann-Krause dynamics, which provides a quite flexible approach to study opinion formation phenomena.

\subsection{Laplacian dynamics on temporal graphs}

A commonly defined collective behavior in Newtonian particle systems is asymptotic flocking. This topic has been extensively studied in the many literature. The definition is as follows:
\begin{definition}
    An interacting many-body system $\{(\bx_i,\bv_i)\}_{i\in[N]}$ exhibits asymptotic flocking if and only if it satisfies the following two conditions:
    \begin{enumerate}
        \item {\bf (Group formation)} The relative distance between particles is uniformly bounded:
        \begin{align*}
            \sup_{0\le t<\infty}\max_{i,j\in[N]}\|\bx_i(t)-\bx_j(t)\|<\infty.
        \end{align*}
        \item {\bf (Velocity alignment)} The relative velocity fluctuations tend to zero as time goes to infinity:
        \begin{align*}
            \lim_{t\to\infty}\max_{i,j\in[N]}\|\bv_i(t)-\bv_j(t)\|=0.
        \end{align*}
    \end{enumerate}
\end{definition}
As discussed in the model derivation part (\Cref{model_deri}), one possibility to analyze our dynamics is to interpret the kernel terms in the adaptive CS model as the edges and weights of a graph. Hence, the model can be viewed as a Laplacian Dynamics or Linear Consensus Model (see \Cref{sec_prelim}) on a temporal (or time-varying) graph and the direction of our analysis follows this interpretation in this work. In the mathematical description, there are paratactic Laplacian dynamics for each dimension $k\in[n]$ in the form of
\begin{align*}
    \dot\bv_{(k)}=L_{\mathcal G(t)}\bv_{(k)},
\end{align*}
where $L_{\mathcal G(t)}$ is the Laplacian matrix of a temporal graph $\mathcal G(t)=([N],\mathcal E(t),\mathcal W(t))$ with the notation
\begin{align*}
    \bv_i=(v_{i1},\cdots,v_{in})^\mathsf{T}\in\mathbb R^n,\quad\bv_{(k)}:=(v_{1k},\cdots, v_{Nk})^\mathsf{T}\in\mathbb R^N.
\end{align*}
\indent To investigate how well the adaptive rule in \eqref{main_model} operates as intended, this paper studies the singular limit as $\varepsilon\to\infty$. By Fenichel's Theorem \cite{Fenichel1979}, this analysis reveals the asymptotic dynamics of the adaptive CS model for sufficiently large $\varepsilon$. Our main novel technique is to reinterpret the nonlinear system as a Laplacian dynamics on a time-varying graph, and to investigate its asymptotic behavior by analyzing the resulting linear system. The key idea is to exploit structural features of the time-varying graph, induced by the evolution of positions and velocities in the adaptive CS model, in order to analyze the corresponding Laplacian dynamics. In short, we extract useful information from each system, combine it, and derive the final result. This methodology is not limited to the present model; it can also be applied to a broader class of systems exhibiting a similar structure, where the effective slow node dynamics can be studied using linear non-autonomous systems.

We focus on the dynamics of one-layer CS system under the assumption that the individual coupling strengths fully reflect the velocity similarity at each time. In this case, the individual coupling strengths are assumed to satisfy $\dot\kappa_{ij}=0$, which leads to symmetry of the coupling strengths. From now on, to avoid any confusion, we use the notation $a_{ij}$ instead of $\kappa_{ij}$, i.e.,
\begin{align*}
    \dot\bx_i=\bv_i,\quad\dot\bv_i=\frac{\kappa}{\mathcal I_i}\sum_{j=1}^N\psi_{ij}a_{ij}(\bv_j-\bv_i),\quad\mbox{where}\quad a_{ij}=\frac{\bv_i\cdot\bv_j}{|\bv_i||\bv_j|}.
\end{align*}
The adjacency matrix of the graph in the singular adaptive CS model is represented by
\begin{align*}
    \left(\left[\frac{\kappa\psi_{ij}(t)a_{ij}(t)}{\mathcal I_i(t)}\right]\right)_{ij}
\end{align*}
and has the following properties:
\begin{enumerate}[label=(\roman*),itemsep=0em]
    \item $\psi_{ij}$ indicates the existence of an {\it undirected} edge between vertices $i$ and $j$,
    \item for a fixed vertex $i$, the velocity similarity $a_{ij}$, which is symmetric, is normalized according to the number of (directed connected) neighbors.
\end{enumerate}
\noindent Accordingly, if $(\psi_{ij})_{ij}$ is non-trivial, the graph is asymmetric in many cases because of normalizing factor $\mathcal I_i$. Under the assumption of edge set changes at most countable times (see \Cref{subsec_existence}), we aim to explore three cases with a more detailed mathematical analysis:

\vspace{.3cm}

\noindent$\blacktriangleright$ \textit{$(\psi_{ij}(t))_{ij}$ is an all-ones matrix:}

\noindent In this scenario, every normalizer is same $\mathcal I_i=N$, hence the graph remains symmetric. If all weights uniformly have a positive lower bound, the Laplacian dynamics on the temporal graph will always converge to the average consensus at an exponential rate, i.e., all nodes converge to the average value of their initial opinions as time tends to infinity. The Fiedler value approach (or spectral analysis) is used for the analysis of this effect; see \Cref{thm_complete}.

\vspace{.3cm}

\noindent$\blacktriangleright$ \textit{$(\psi_{ij}(t))_{ij}$ has neighbor-connectivity (See \Cref{def_neighbor_conn}):}

\noindent In this case, there is no longer any room for discussion of symmetry. Nonetheless, if all weights uniformly have a positive lower bound, the neighbor-connectivity guarantees a negative dichotomy spectrum \cite{Fink2006} in the evolution of relative velocities. Consequently, the convergence of relative velocities can be reached at exponential rate; see Theorem~\Cref{a_priori_thm}.

\vspace{.3cm}

\noindent$\blacktriangleright$ \textit{The graph is strongly connected and ``almost" symmetric:}

\noindent The asymmetry arises from different normalizing factors $\mathcal I_i$. To handle the asymmetric part in graphs, we assume that the factors have relatively close lower and upper bounds. Using the incidence structure of our graphs, one can provide a rather technically challenging framework to still obtain fast convergence to asymptotic consensus. (See \Cref{asymm_thm})

\vspace{.3cm}

The rest of this paper is organized as follows. In \Cref{sec_prelim}, we mention the existence and uniqueness of solutions to the adaptive CS model, which is a system of discontinuous ordinary differential equations. And we recall the basic graph theory setup and relevant linear algebra concepts. In \Cref{sec_thm}, we present three main results of the singular adaptive CS model or Laplacian dynamics on three different types of temporal graphs. \Cref{sec_ex} provides an analysis of a two-particle toy model for flocking and dispersion, and offers examples of temporal graphs discussed in the main results of \Cref{sec_thm}. In \Cref{sec_num}, we present numerical simulations that support the main theorems and explore fully adaptive CS dynamics numerically beyond the singular limit analysis. Finally, \Cref{sec_concl} is devoted to a brief summary of the main results and some open problems.

\section{Preliminaries}\label{sec_prelim}
In this section, we discuss the existence and uniqueness of solution to adaptive CS model \eqref{main_model} in the sense of Carath\'eodory solution, which contains discontinuous ODEs, and recall the basic concepts in graph theory and linear algebra.
\subsection{Existence and uniqueness of solution}\label{subsec_existence}

The velocity vectors in the adaptive CS model \eqref{main_model} evolve according to a discontinuous ODE. The existence and uniqueness of solutions to discontinuous ODEs is a mathematically challenge, which has given rise to several notions of generalized solutions, including Carath\'eodory, Filippov, and Krasovskii solutions. 

Most discontinuous ODEs do not admit classical solutions if the solution trajectory passes through points of discontinuity in the vector field. The Carath\'eodory solution offers the most immediate generalization of the classical notion. Roughly speaking, a Carath\'eodory solution exists if the discontinuities occur at most at countably many points; however, uniqueness is generally not ensured. 

In the case of the continuous-time Hegselmann–Krause model, Blondel, Hendrickx, and Tsitsiklis \cite{BHT2010} proved both existence and uniqueness of Carath\'eodory solutions via detailed estimates. They explicitly constructed an almost sure set of initial conditions, and showed that if the dynamics start from a point in this set, then the solution necessarily transitions to the neighboring vector field whenever a discontinuity occurs. This ensures uniqueness and allows the dynamics to be analyzed in terms of ODE expressions except on a set of times of measure zero. 

Ceragioli and Frasca \cite{CF2012} established the existence of Krasovskii solutions for the same model and proposed a framework to analyze its asymptotic behavior within this more general setting. Their results also involved relaxed conditions on the initial set compared to those in \cite{BHT2010}. 

One decade after these works, Piccoli and Rossi \cite{PR2021} considered a generalization of the same model with nonlinear interaction kernels, and rigorously clarified the relationship among classical, Carath\'eodory, Filippov, and Krasovskii solutions. They proved existence and uniqueness results in each setting. However, these results rely heavily on the model’s structure and therefore do not permit straightforward generalization.

It is also worth noting that, half a century ago, H\'ajek \cite{Hajek1979} showed that if the discontinuous vector field is locally bounded, then a local Krasovskii solution exists for every initial condition. This implies, for example, that in the Kuramoto-type model with the connection kernel and normalizing factor we propose, existence is automatically guaranteed, in contrast to the CS model.

In this work, we do not address the existence or uniqueness of solutions for our model as it forms really a technically very separate issue from the analysis of the dynamics. Instead, all of our analysis proceeds under the assumption that there exists a set of initial conditions for which the adaptive CS model admits a unique global Carath\'eodory solution, and we restrict attention to initial data within this set. This assumption implies that, in the Laplacian dynamics considered in our main results, the network's edge set $\mathcal E(t)$ changes at most countably many times over time. Hence, the Carath\'eodory condition holds so that the existence and uniqueness of a Carath\'eodory solution are ensured. A rigorous investigation of the existence and measure-theoretic properties of the proper initial data set for the adaptive CS model will be the subject of future separate work.

\subsection{Basic theory}
\label{sec_prelim}

In this subsection, we recall the setup in graphs, fundamental properties of Laplacian dynamics. A graph $\mathcal G=(\mathcal V,\mathcal E,\mathcal W)$ consists of a set of nodes $\mathcal V=\{\sx_1,\cdots,\sx_N\}$, a set of edges $\mathcal E\subseteq\mathcal V\times\mathcal V$, and a weight function $\mathcal W:\mathcal W\to\mathbb R$ that assigns a real-valued weight to each edge. The weight $w_{ij}=\mathcal W(\sx_i,\sx_j)$ quantifies the strength or influence of the connection between nodes $\sx_i$ and $\sx_j$. Positive weights indicate cooperative interactions, while negative weights represent antagonistic interactions. Edges with zero weight are considered absent.

We consider both undirected and directed graphs, allowing for weighted and signed edges. In an \textit{undirected graph}, each edge $(\sx_i,\sx_j)\in\mathcal E$ implies $(\sx_j,\sx_i)\in\mathcal E$, and the weights satisfies $w_{ij}=w_{ji}$. In a \textit{directed graph}, each edges have a specified direction, meaning $(\sx_i,\sx_j)\in\mathcal E$ does not necessarily imply $(\sx_j,\sx_i)\in\mathcal E$ and asymmtry of weights. A graph can also be represented by an \textit{adjacency matrix} $A$, where $A_{ij}=w_{ij}$ if $(\sx_i,\sx_j)\in\mathcal E$, and $A_{ij}=0$ otherwise. For undirected graphs, $A$ is symmetric while for directed graphs, $A$ is generally asymmetric. 

The degree of a node quantifies the total weight of its connections. For undirected graphs, the \textit{degree} of node $i$ is given by
\begin{align*}
    \mbox{deg}(\sx_i)=\sum_{j=1}^NA_{ij}.
\end{align*}
For directed graphs, we use the \textit{out-degree} defined by
\begin{align*}
    \mbox{deg}_{\text{out}}(\sx_i)=\sum_{j=1}^NA_{ij}.
\end{align*}
The degree matrices for undirected and directed graphs are defined as
\[D=\mbox{diag}(\mbox{deg}(\sx_1),\cdots,\mbox{deg}(\sx_N))\quad\mbox{and}\quad D_{\mbox{out}}=\mbox{diag}(\mbox{deg}_{\text{out}}(\sx_1),\cdots,\mbox{deg}_{\text{out}}(\sx_N)),\]
respectively. To simplify the notation, we denote both matrices by $D$ and refer to them as degree matrices, as they are identical in their mathematical form. The \textit{Laplacian matrix}, originating from the discrete diffusion operator in finite difference methods, is defined by
\begin{align*}
    L=D-A.
\end{align*}
Laplacian dynamics describes the evolution of node states on a graph, where each node interacts with its neighbors to minimize the state differences over time. In the case of a signed graph, negative weights represent opposing interactions, which reverse the typical effect. The dynamics is governed by
\begin{align*}
    \dot\sX=-L\sX,
\end{align*}
where $\sX=[\sx_1(t),\cdots,\sx_N(t)]^{\mathsf{T}}\in\mathbb R^N$ is the state vector of the nodes, and $N$ denotes the total number of nodes. The asymptotic behavior of the node states is closely tied to the eigenvalues of the Laplacian matrix $L$. On static graphs, extensive studies have been conducted to analyze Laplacian dynamics in this context. The following theorem summarizes key properties of the Laplacian matrix, which are essential for analyzing Laplacian dynamics.

\begin{theorem}\label{thm_Lapla_prop}
    The Laplacian matrix of a graph has smallest eigenvalue zero with the corresponding eigenvector $\mathds{1}=(1,\cdots,1)^\mathsf{T}$ and algebraic multiplicity equal to the number of components of the graph.
\end{theorem}
Next, we provide a useful expression for the analysis of the signs of the Laplacian matrix.
\begin{lemma}\label{lem_vector_weight}
    For any weighted undirected graph $\mathcal G=(\mathcal V,\mathcal E,\mathcal W)$ with adjacency matrix $A$, one has
    \[x^TLx=\sum_{i,j\in\mathcal V}\sqrt{a_{ij}a_{ji}}(x_i-x_j)^2,\quad\forall~x=(x_1,\cdots x_n)^T\in\mathbb R^n.\]
\end{lemma}
\begin{proof}
    To avoid any confusion, we clarify that
    \[\mathcal E\subseteq\{(i,j)~:~i,j\in\mathcal V\mbox{ such that }i>j\}.\]
    Consider the matrix $\nabla$ defined by as follows: for $e=(i,j)\in\mathcal E$,
    \begin{align*}
        \nabla_{e,i}=\sqrt{a_{ij}},\quad\nabla_{e,j}=-\sqrt{a_{ij}}\quad\mbox{and}\quad\nabla_{e,k}=0,~~\forall~k\ne i,j.
    \end{align*}
    For the Laplacian matrix $L$ of $\mathcal G$ and a vector $x=(x_1,\cdots,x_N)^\mathsf{T}\in\mathbb R^N$, the matrix $\nabla$ satisfies
    \begin{align*}
        L=\nabla^\mathsf{T}\nabla,\quad(\nabla x)_{e=(i,j)}=\sqrt{a_{ij}}(x_i-x_j),
    \end{align*}
    where the symmetry of $A$ is used. Hence, we have
    \begin{align*}
        x^\mathsf{T}Lx=\sum_{\substack{i,j\in\mathcal V\\i>j}}a_{ij}(x_i-x_j)^2=\frac{1}{2}\sum_{i,j\in\mathcal V}a_{ij}(x_i-x_j)^2.
    \end{align*}
    This ends the proof.
\end{proof}
For more information about Laplacian matrix, we refer to \cite{GR2001}.
\section{Asymptotic consensus on temporal graph}
\label{sec_thm}

In this section, we analyze the asymptotic behavior of the Cauchy problem for the fast-adaptation singular limit of the adaptive CS model
\begin{align}\label{singular_adapCS}
    \begin{aligned}
    \begin{dcases}
        \dot\bx_i=\bv_i,\quad\dot\bv_i=\frac{\kappa}{\mathcal I_i}\sum_{j=1}^N\psi_{ij}a_{ij}(\bv_j-\bv_i),\\
        \bx_i(0)=\bx_i^0,~~\bv_i(0)=\bv_i^0,\quad\forall~i\in[N],
    \end{dcases}
    \end{aligned}
\end{align}
where $a_{ij}$ is defined by
\[a_{ij}(t):=\frac{\bv_i(t)\cdot\bv_j(t)}{|\bv_i(t)||\bv_j(t)|},\quad\forall~i,j\in[N].\]
In the following three subsections, we present the main results concerning Laplacian dynamics on temporal graphs
\begin{align}\label{CS_graph}
    \left(\frac{\kappa}{\mathcal I_i}\psi_{ij}a_{ij}\right)_{ij}
\end{align}
and demonstrate their application to the singular adaptive CS model \eqref{singular_adapCS}.

In the first subsection, which deals with the complete network where no discontinuity occurs in the vector field, the Cauchy–Lipschitz theorem guarantees the existence of a unique classical solution for any initial data. However, as discussed in \Cref{subsec_existence}, in the following two subsections we consider the set of initial data generating unique global Carathéodory solutions to \eqref{singular_adapCS}, denoted by $\mathcal D$. The initial data is assumed to belong to $\mathcal D$.
\subsection{Laplacian dynamics on complete networks}\label{subsec_complete}
In this subsection, we consider Laplacian dynamics on a temporal graph $\mathcal G(t)=(\mathcal V,\mathcal E,\mathcal W(t))$, which is a complete network, i.e., $\mathcal E$ is time-independent with
\begin{align*}
    (i,j)\in\mathcal E,\quad\forall~i,j\in[N].
\end{align*}
We restrict our interest to graphs whose weights are symmetric and have a positive uniform lower bound. In other words, for particles $\sx_i\in\mathbb R$, $i\in[N]$, we study
\begin{align*}
    \dot\sx_i=\sum_{j=1}^N\alpha_{ij}(t)(\sx_j-\sx_i),\quad\forall~i\in[N],
\end{align*}
where the temporal network satisfies
\begin{align}\label{weight_bound}
    \alpha_m\le\alpha_{ij}(t),\quad\alpha_{ij}(t)=\alpha_{ji}(t),\quad\forall~t\ge0,~~\forall~i,j\in[N],
\end{align}
for a positive constant $\alpha_m>0$. The above system can be written in a linear form
\begin{align}\label{system1}
    \dot\sX(t)=-L(t)\sX(t),
\end{align}
where $\sX=(\sx_1,\cdots,\sx_N)^\mathsf{T}\in\mathbb R^N$ and $L(t)$ denotes the Laplacian matrix of $\mathcal G(t)$. Note that the one-dimensional study does not lose generality in dimension.

The main idea to deal with this temporal network in dynamics is to bound the eigenvalues of Laplacian matrices by comparing those of other static networks. In the following lemma, we describe the spectrum relation between two static graphs whose adjacency matrices satisfy an element-wise matrix inequality.
\begin{lemma}\label{eigval_comparison}
    Let $\mathcal G_i=(\mathcal V,\mathcal E_i,\mathcal W_i)$ for $i=P,Q$ be two weighted undirected static graphs with adjacency matrices $A_P=(p_{ij})$ and $A_Q=(q_{ij})$, respectively. If
    \begin{align*}
        0\le p_{ij}\le q_{ij},\quad\forall~i,j\in[N],
    \end{align*}
    then the eigenvalues $\{\lambda_k^i\}_k\in[N]$ for $i=P,Q$, of their respective Laplacian matrices, arranged in ascending order, satisfy
    \begin{align*}
        \lambda_k^P\le\lambda_k^Q,\quad\forall~k\in[N].
    \end{align*}
\end{lemma}
\begin{proof}
    Since $A_P$ and $A_Q$ are symmetric, by \Cref{lem_vector_weight}, we have
    \begin{align*}
        x^\mathsf{T}L_Px=\sum_{i,j\in\mathcal V}p_{ij}(x_i-x_j)^2,\quad\mbox{similarly,}~~x^\mathsf{T}L_Qx=\sum_{i,j\in\mathcal V}q_{ij}(x_i-x_j)^2.
    \end{align*}
    From element-wise matrix inequality, they satisfy
    \begin{align*}
        x^\mathsf{T}L_Px\le x^\mathsf{T}L_Qx,\quad\forall~x\in\mathbb R^N.
    \end{align*}
    Now, we use Courant-Fishcer Theorem \cite{HJ1985} to derive the desired result.    
\end{proof}
\begin{theorem}\label{thm_complete}
    The Laplacian dynamics \eqref{system1} with the assumption \eqref{weight_bound} exhibits asymptotic average consensus. In particular, the time-averaged value of the particles, denoted by $\bar\sx$, is constant, and the convergence is exponential,
    \begin{align*}
        \sum_{i=1}^N(\sx_i(t)-\bar\sx)^2\le e^{-2\alpha_mNt}\sum_{i=1}^N(\sx_i(0)-\bar\sx)^2,\quad\forall~t\ge0,~~\forall~i\in[N].
    \end{align*}
\end{theorem}
\begin{proof}
    \noindent$\bullet$ {\bf (Constant average $\bar\sx$)} Denote by $\bar\sx$ the average of all particles, that is, $\bar\sx=\mathds{1}^\mathsf{T}\sX$ for one vector $\mathds{1}=(1,\cdots,1)^\mathsf{T}\in\mathbb R^N$. Its dynamics follows
    \begin{align*}
        \dot{\bar\sx}&=\mathds{1}^\mathsf{T}\dot\sX=-\mathds{1}^\mathsf{T}L(t)X(t)=-L(t)\mathds{1}X(t)=0,
    \end{align*}
    where we used symmetry of $L$ and $L\mathds{1}=\mathbf{0}$. Hence, we have $\bar\sx(t)\equiv\bar\sx(0)$.

    \vspace{.2cm}

    \noindent$\bullet$ {\bf (Asymptotic consensus)} Now, we consider the fluctuation vector $\tilde\sX\in\mathbb R^N$ defined by
    \[\tilde\sX:=(\sx_1-\bar\sx,\cdots,\sx_N-\bar\sx)^\mathsf{T},\]
    and rewrite the system \eqref{system1} as
    \begin{align}\label{fluc_sys}
        \dot{\tilde\sX}(t)=-L(t)\tilde\sX(t).
    \end{align}
    Then, we have a Lyapunov function
    \[F(t):=\frac{1}{2}\|\tilde\sX(t)\|^2\]
    for \eqref{fluc_sys}:
    \begin{align*}
        \frac{d}{dt}F(t)&=-\tilde\sX L(t)\tilde\sX\le-\min_{\substack{\sY\ne\mathbf{0}\\\mathds{1}^\mathsf{T}\sY=0}}\frac{\sY^\mathsf{T}L(t)\sY}{\|\sY\|^2}\cdot\|\tilde\sX\|^2=-\lambda_2(L(t))\|\tilde\sX\|^2\le0,
    \end{align*}
    where we used $\mathds{1}^\mathsf{T}\tilde\sX=0$ and the fact from \Cref{thm_Lapla_prop} and Courant-Fischer theorem that
    \begin{align*}
        \min_{\substack{\bx\ne0\\\mathds{1}^\mathsf{T}\bx=0}}\frac{\bx^\mathsf{T}L\bx}{\|\bx\|^2}=\lambda_2(L)
    \end{align*}
    for the spectrum $(\lambda_i)$ in ascending order. Here, we can use \Cref{eigval_comparison} to compare $\lambda_2(L(t))$ with a graph which has adjacency matrix $\alpha_mJ_N$, where $J_N\in\mathbb R^{N\times N}$ denotes the all-ones matrix. A direct computation of the derivative gives
    \begin{align*}
        \frac{d}{dt}\|\tilde \sX\|\le-\alpha_mN\|\tilde\sX\|,
    \end{align*}
    which implies the desired result.
\end{proof}
We now demonstrate the application of \Cref{thm_complete} in our model to derive asymptotic average consensus. First, we show the boundedness of network weights.
\begin{lemma}\label{a_bound}
    Suppose that initial data $\{(\bx_i^0,\bv_i^0)\}_{i\in[N]}$ satisfy
    \[a_m:=\min_{i,j\in[N]}\frac{\bv_i^0\cdot\bv_j^0}{|\bv_i^0||\bv_j^0|}>0.\]
    Then, the solution $\{(\bx_i,\bv_i)\}_{i\in[N]}$ to \eqref{singular_adapCS} satisfies
    \[\frac{\bv_i\cdot\bv_j}{|\bv_i||\bv_j|}\ge a_m,\quad\forall~t\ge0.\]
\end{lemma}
\begin{proof}    
    Recall that the geometric interpretation of dot product gives
    \begin{align}\label{def_a}
        a_{ij}=\frac{\bv_i\cdot\bv_j}{|\bv_i||\bv_j|}=\cos\left(\measuredangle(\bv_i,\bv_j)\right),
    \end{align}    
    where $\measuredangle(\bv_i,\bv_j)$ denotes the angle between $\bv_i$ and $\bv_j$. Since we have the positivity of $a_{ij}$ for any pair of initial data, the largest angle at $t=0$ is less than $\pi/2$. Thus, the following claim guarantees the monotone increase of $\min_{i,j\in[N]}a_{ij}(t)$ in time.
    \begin{center}
        {\it Claim:} $\mbox{conv}\left(\{\bv_i(s)\}_{i\in[N]}\right)\subseteq\mbox{conv}\left(\{\bv_i(t)\}_{i\in[N]}\right)$,$\quad\forall~s\ge t$.
    \end{center}
    Here, $\mbox{conv}(\cdot)$ denotes the convex hull of a vector set.\newline
    \noindent{\it Proof of claim:} Suppose not. Then, the time point $t_0$ defined by
    \[t_0:=\inf\left\{s\ge t~:~\mbox{conv}\left(\{\bv_i(s)\}_{i\in[N]}\right)\not\subseteq\mbox{conv}\left(\{\bv_i(t)\}_{i\in[N]}\right)\right\}\]
    is finite. Since a convex hull is a closed set, there exists a constant $0<h\ll1$ such that
    \begin{align}\label{conv_eq}
        \begin{aligned}
            \mbox{conv}\left(\{\bv_i(t_0)\}_{i\in[N]}\right)&\subseteq\mbox{conv}\left(\{\bv_i(t)\}_{i\in[N]}\right),\\
            \mbox{conv}\left(\{\bv_i(t_0+\tilde h)\}_{i\in[N]}\right)&\not\subseteq\mbox{conv}\left(\{\bv_i(t)\}_{i\in[N]}\right),\quad\forall~0<\tilde h<h.
        \end{aligned}
    \end{align}
    The relations \eqref{conv_eq} implies that there exist an index $k\in[N]$ and a constant $\delta>0$ such that
    \begin{align}\label{conv_eq2}
        \bv_k(t_0)+\tilde\delta\dot\bv_k(t_0)\notin\mbox{conv}\left(\{\bv_i(t_0)\}_{i\in[N]}\right),\quad\forall~0<\tilde\delta<\delta.
    \end{align}
    On the other hand, note that
    \begin{align*}
        a_{ij}(t)\in[0,1]\quad\mbox{and}\quad0\le\frac{1}{N}\sum_{j=1}^N\psi_{ij}(t)a_{ij}(t)\le1.
    \end{align*}
    Thus, from the governing equation \eqref{singular_adapCS} and the definition of a convex hull, one can easily see that there exists a constant $0<\varepsilon\ll1$ such that
    \begin{align*}
        \bv_j(t)+\tilde\varepsilon\dot\bv_j(t)\in\mbox{conv}\left(\{\bv_i(t)\}_{i\in[N]}\right),\quad\forall~j\in[N],~~\forall~0<\tilde\varepsilon<\varepsilon,
    \end{align*}
    at any fixed $t>0$. This derives a contradiction to \eqref{conv_eq2}, which ends the proof.
\end{proof}
\begin{figure}
    \centering
    \includegraphics[width=0.6\linewidth]{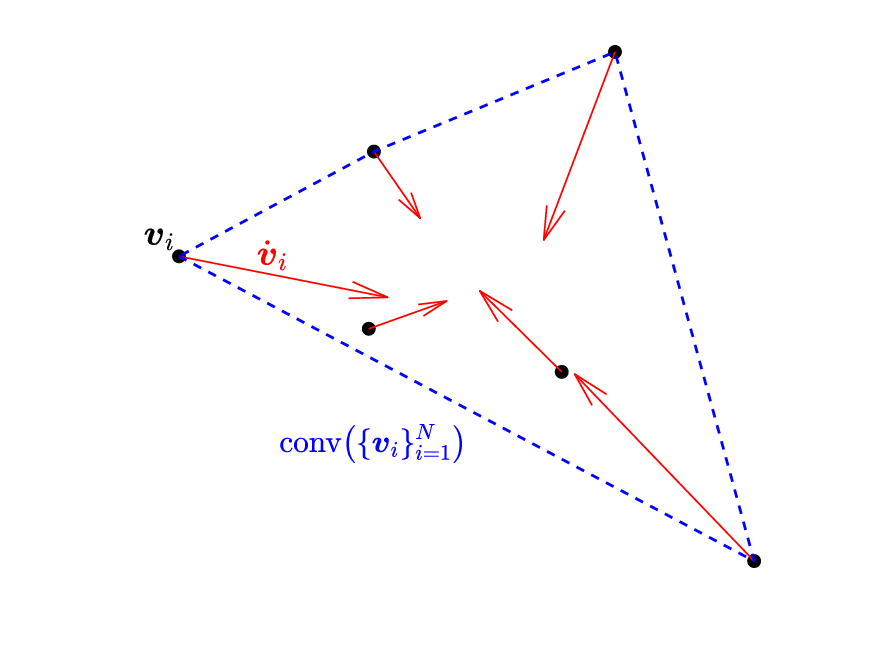}
    \caption{Sketch of proof in \Cref{a_bound}}
    \label{fig:a_bound}
\end{figure}
We define the average velocity and fluctuation ensemble by
\begin{align*}
    \bv_c(t):=\frac{1}{N}\sum_{i=1}^N\bv_i(t)\in\mathbb R^n,\quad\tilde\bv:=(\bv_1-\bv_c,\cdots,\bv_N-\bv_c)\in\mathbb R^{n\times N}.
\end{align*}
\begin{corollary}\label{full_connected}
    Suppose that initial data $\{(\bx_i^0,\bv_i^0)\}_{i\in[N]}$ satisfy
    \begin{align*}
        \max_{i,j\in[N]}\|\bx_i^0-\bx_j^0\|\le d_0,\quad a_m:=\min_{i,j\in[N]}\frac{\bv_i^0\cdot\bv_j^0}{|\bv_i^0||\bv_j^0|}>0,
    \end{align*}
    for a positive constant $d_0<d-2\|\tilde{\bv}(0)\|_F/a_m$, where $\|\cdot\|_F$ denotes the Frobenius norm. Then, the solution $\{(\bx_i,\bv_i)\}_{i\in[N]}$ to \eqref{singular_adapCS} exhibits asymptotic flocking.
\end{corollary}
\begin{proof}
    First, we show that all particles are within a detection range each other, i.e.,
    \begin{align*}
        \psi_{ij}(t)=1,\quad\forall~i,j\in[N],~~\forall~t\ge0.
    \end{align*}
    Set a constant
    \begin{align*}
        \mathcal T:=\sup\{t>0~:~\psi_{ij}(t)=1,\quad\forall~i,j\in[N],~~\forall~s\in[0,t)\}.
    \end{align*}
    By the initial condition, $\mathcal T>0$ holds. On the time interval $t\in[0,\mathcal T)$, we may write the system \eqref{singular_adapCS} as
    \begin{align}\label{HK_2}
        \begin{aligned}
            \begin{dcases}
                \dot\bx_i=\bv_i,\\
                \dot\bv_i=\sum_{j=1}^N\frac{a_{ij}}{N}(\bv_j-\bv_i),
            \end{dcases}
        \end{aligned}
    \end{align}
    where $a_{ij}$ is defined in \eqref{def_a}. We claim that the completely connected structure remains over the time.
    \begin{center}
        {\it Claim:} $\mathcal T=\infty$.
    \end{center}
    \noindent{\it Proof of claim:} Suppose not. Then, there exists an index pair $i^*,j^*\in[N]$ such that
    \begin{align}\label{claim_assertion}
        \|\bx_{i^*}(\mathcal T)-\bx_{j^*}(\mathcal T)\|=d.
    \end{align}
    By \Cref{a_bound}, we have
    \[a_m\le a_{ij}(t)\le1,\quad\forall~t\ge0,~~\forall~i,j\in[N].\]
    Observe that the system \eqref{HK_2} can be rewritten as Laplacian-type dynamics if we decompose the components of the velocity vector $\bv_i\in\mathbb R^d$ and observe that each component satisfies Laplacian dynamics of the form \eqref{system1}. This leads to $d$ Laplacian-type differential equations. Therefore, we can apply \Cref{thm_complete} to derive
    \begin{align*}
        \|\tilde\bv(t)\|_F\le e^{-a_mt}\|\tilde\bv(0)\|_F,\quad\forall~t\in[0,\mathcal T).
    \end{align*}
    This implies the individual convergence rate to the average vector, i.e.,
    \begin{align}\label{dec_rate}
        \|\bv_i-\bv_c\|\le e^{-a_mt}\|\tilde\bv(0)\|_F,\quad\quad\forall~i\in[N],~~\forall~t\in[0,\mathcal T).
    \end{align}
    We combine \eqref{dec_rate} and the relation
    \begin{align*}
        \frac{d}{dt}\|\bx_i-\bx_j\|\le\|\bv_i-\bv_j\|
    \end{align*}
    to derive
    \begin{align*}
        \|\bx_i(t)-\bx_j(t)\|\le d_0+2\|\tilde\bv(0)\|_F\int_0^te^{-\kappa_ms}~\textnormal{d}s<d,\quad\forall~t\in[0,\mathcal T),
    \end{align*}
    for any $i,j\in[N]$. By the continuity of relative distances, we have
    \begin{align*}
        \|\bx_i(\mathcal T)-\bx_j(\mathcal T)\|<d,\quad\forall~i,j\in[N],
    \end{align*}
    which is a contradiction to \eqref{claim_assertion}. This ends the proof of claim.

    \vspace{.2cm}    

    \noindent Now, we can show the velocity alignment for asymptotic flocking. By the claim, the inequality \eqref{dec_rate} holds for all $t\ge0$. Therefore, we have the fast convergence of velocities,
    \begin{align*}
        \lim_{t\to\infty}\|\bv_i(t)-\bv_c(0)\|=0,\quad\forall~i\in[N],
    \end{align*}
    which finishes the proof.
\end{proof}
\subsection{Laplacian dynamics on neighbor-connected networks}
In this subsection, we study Laplacian dynamics on temporal digraphs that have a certain type of connectivity. To be more specific, we consider digraphs where any two vertices have distance at most 2. In the following definition, we clarify the details.
\begin{definition}\label{def_neighbor_conn}
    A directed graph $\mathcal G=(\mathcal V,\mathcal E)$ is said to be neighbor-connected if for any two distinct vertices $i,j\in\mathcal V$, at least one of the following conditions holds:\newline
    (i) $i$ and $j$ are directly connected by edges, i.e., $(i,j),(j,i)\in\mathcal E$,\newline
    (ii) There exists a vertex $k\in\mathcal V$ such that $(i,k),(k,i),(j,k),(k,j)\in\mathcal E$.
\end{definition}
From now on, we consider the Laplacian dynamics on a temporal graph $\mathcal G(t)=(\mathcal V,\mathcal E(t),\mathcal W(t))$ satisfying the following three conditions $(\mathcal A)$:

\vspace{.1cm}

\begin{enumerate}[label={($\mathcal A$\arabic*)}]
    \item for each $t\ge0$, the graph $\mathcal G(t)$ is neighbor-connected,
    \item weights $\mathcal W(t)$ have a positive uniform lower bounds, i.e.,
    \begin{align*}
        \exists~w_m>0~~\mbox{such that}~~w_m\le w_{ij}(t),~~\forall~t\ge0,~~\forall~(i,j)\in\mathcal E(t),
    \end{align*}
    \item the edge set $\mathcal E(t)$ changes at most countable times.
\end{enumerate}

\vspace{.1cm}

In the mathematical description, we have the following governing equations: 
\begin{align}\label{system2}
    \dot\sx_i=\sum_{j=1}^N\beta_{ij}(t)(\sx_j-\sx_i),\quad\forall~i\in[N],
\end{align}
where $(\beta_{ij}(t))$ is the adjacency matrix of $\mathcal G(t)$. In this case, we will work directly with multi-dimensional vectors $\sx_i\in\mathbb R^n$. In the following theorem, we show the asymptotic consensus of \eqref{system2}.
\begin{theorem}\label{a_priori_thm}
    The Laplacian dynamics \eqref{system2} with the assumptions $(\mathcal A)$ exhibits asymptotic consensus. Moreover, the convergence rate is exponential, i.e., there exists $C>0$ such that
    \begin{align*}
        \|\sx_i(t)-\sx_j(t)\|\le Ce^{-w_mt},\quad\forall~t\ge0,~~\forall~i,j\in[N].
    \end{align*}
\end{theorem}
\begin{proof}
    In this proof, we derive estimates based on the Carath\'eodory solution. Its existence and uniqueness are guaranteed by ($\mathcal A3)$. The key strategy is to show the convergence of differences between particles using the dichotomy spectrum \cite{JKL2024}. Let $\eta:\mathbb R\to\mathbb R_{\ge0}$ be the positive part function defined by
    \begin{align*}
        \eta(x)=
        \begin{cases}
            x,\quad&\mbox{if}~~x\ge0,\\
            0,&\mbox{otherwise,}
        \end{cases}
    \end{align*}
    and consider the synchronization errors given by the vector of $N(N-1)/2$ components
    \begin{align*}
        \xi(t)=(\xi_{ij}(t))^T_{i,j\in[N],i<j}\quad\mbox{where}\quad\xi_{ij}(t)=\|\sx_i(t)-\sx_j(t)\|^2.
    \end{align*}
    From \eqref{system2}, one can derive
    \begin{align*}
        \frac{d}{dt}\|\sx_i-\sx_j\|^2&=-\left(2(a_{ij}(t)+a_{ji}(t))+\sum_{m\ne i,j}(a_{im}(t)+a_{jm}(t))\right)\|\sx_i-\sx_j\|^2\\
        &\hspace{.4cm}-\sum_{m\ne i,j}(a_{im}(t)-a_{jm}(t))\left(\|\sx_i-\sx_m\|^2-\|\sx_j-\sx_m\|^2\right).
    \end{align*}
    Hence, by using the function $\eta$, we can get a differential inequality
    \begin{align*}
        \frac{d}{dt}\xi(t)\le E(t)\xi(t),
    \end{align*}
    where the time-dependent matrix $E(t)$ is defined rowwise as follows: let us label each of the $N(N-1)/2$ rows by $e^{ij}\in\mathbb R^{N(N-1)/2}$, where $i,j\in[N]$ with $i<j$, i.e.,
    \begin{align*}
        e^{ij}(t)=(e^{ij}_{lk}(t))_{l,k\in[N], l<k}.
    \end{align*}
    For fixed $i,j$, the elements satisfy
    \begin{align*}
        e_{lk}^{ij}(t)=
        \begin{dcases}
            -2(a_{ij}(t)+a_{ji}(t))-\sum_{\tilde l\ne i,j}(a_{i\tilde l}(t)+a_{j\tilde l}(t)),~~&\mbox{if}~~l=i~\mbox{and}~k=j,\\
            \eta\left(a_{jk}(t)-a_{ik}(t)\right),&\mbox{if}~~l=i~\mbox{and}~k\ne j,\\
            \eta\left(a_{jl}(t)-a_{il}(t)\right),&\mbox{if}~~k=i,\\
            \eta\left(a_{ik}(t)-a_{jk}(t)\right),&\mbox{if}~~l=j,\\
            \eta\left(a_{il}(t)-a_{jl}(t)\right),&\mbox{if}~~l\ne i~\mbox{and}~k=j,\\
            0,&\mbox{otherwise.}
        \end{dcases}
    \end{align*}
    Below, we use a decomposed notation
    \begin{align*}
        a_{ij}(t)=\psi_{ij}(t)w_{ij}(t)    
    \end{align*}
    where 
    \begin{align*}
        \psi_{ij}(t)=
        \begin{cases}
            1,\quad\mbox{if}~~(i,j)\in\mathcal E(t),\\
            0\quad\mbox{otherwise},
        \end{cases}
    \end{align*}
    and $w_{ij}(t)$ is the weight of the edge $(i,j)$. If there is no corresponding edge in $\mathcal W(t)$, any number can be assigned to $w_{ij}(t)$. To simplify the notation, we omit the time dependency.

    Note that the real-valued matrix $E$ has nonnegative entries except for strictly negative diagonal entries. Besides, one can see that $E$ is row-dominant matrix:
    \begin{align*}
        |e^{ij}_{ij}|-\sum_{\substack{l,k\ne i,j\\ l<k}}|e^{ij}_{lk}|&=2(\psi_{ij}w_{ij}+\psi_{ji}w_{ji})+\sum_{\substack{\tilde l\in[N]\\\tilde l\ne i,j}}(\psi_{i\tilde l}w_{i\tilde l}+\psi_{j\tilde l}w_{j\tilde l})\\
        &\hspace{.4cm}-\sum_{\substack{\tilde l\in[N]\\\tilde l\ne i,j}}\Big(\eta(\psi_{j\tilde l}w_{j\tilde l}-\psi_{i\tilde l}w_{i\tilde l})+\eta(\psi_{i\tilde l}w_{i\tilde l}-\psi_{j\tilde l}w_{j\tilde l})\Big)\\
        &=2(\psi_{ij}w_{ij}+\psi_{ji}w_{ji})+\sum_{\substack{l\in[N]\\l\ne i,j}}(\psi_{il}w_{l}+\psi_{jl}w_{jl})-\sum_{\substack{l\in[N]\\l\ne i,j}}|\psi_{il}w_{il}-\psi_{jl}w_{jl}|\\
        &=2(\psi_{ij}w_{ij}+\psi_{ji}w_{ji})+2\left(\sum_{b_l<0}\psi_{il}w_{il}+\sum_{b_l\ge0}\psi_{jl}w_{jl}\right),
    \end{align*}
    where $b_l$ for $l\in[N]\setminus\{i,j\}$ is defined by
    \begin{align*}
        b_l:=\psi_{il}w_{il}-\psi_{jl}w_{jl}.
    \end{align*}
    Thus, we have
    \begin{align}\label{dich_estimate}
        |e^{ij}_{ij}|-\sum_{\substack{l,k\ne i,j\\ l<k}}|e^{ij}_{lk}|\ge2w_m\left(\psi_{ij}+\psi_{ji}+\sum_{b_l<0}\psi_{il}+\sum_{b_l\ge0}\psi_{jl}\right)\ge 2w_m,
    \end{align}
    where we used the neighbor-connected property in the last inequality.

    In consequence, one can see that a linear ODE $\dot u=E(t)u$ has a dichotomy spectrum contained in $(-\infty,0)$ for all $t\ge0$. Similarly as in \cite{JKL2024}, by using a comparison theorem \cite{OO1960} and an exponential dichotomy theorem \cite{Fink2006}, one can derive
    \begin{align*}
        \|\xi_{ij}(t)\|=\|\sx_i(t)-\sx_j(t)\|^2\lesssim e^{-2w_mt},\quad\forall~i<j.
    \end{align*}
    This finishes the proof.
\end{proof}

Now, we are ready to provide asymptotic flocking in \eqref{singular_adapCS} given suitable assumptions hold.

\begin{corollary}\label{cor_thm2}
    Suppose that initial data $\{\bx_i^0,\bv_i^0)\}_{i\in[N]}\in\mathcal D$ satisfy
    \begin{align*}
        w_m:=\min_{i,j\in[N]}\frac{\bv_i^0\cdot\bv_j^0}{|\bv_i||\bv_j|}>0
    \end{align*}
    and assume that the graph \eqref{CS_graph} is neighbor-connected for each $t\ge0$. Then, the solution $\{(\bx_i,\bv_i)\}_{i\in[N]}$ to \eqref{singular_adapCS} exhibits asymptotic flocking.
\end{corollary}
\begin{proof}
     By the initial data and \Cref{a_bound}, $w_m$ is a lower bound of $(w_{ij}(t))_{i,j,t}$. From \Cref{a_priori_thm}, there exists $C>0$ such that
     \begin{align*}
         \|\bv_i(t)-\bv_j(t)\|\le Ce^{-2w_mt},\quad\forall~t\ge0,~~\forall~i,j\in[N].
     \end{align*}
     This exponential decay implies 
     \begin{align*}
         \max_{i,j\in[N]}\sup_{t\ge0}\|\bx_i(t)-\bx_j(t)\|<\infty,
     \end{align*}
     which completes the proof.
\end{proof}
\begin{remark}
    The a priori conditions in \Cref{cor_thm2} are sufficiently realizable and cover a wide range of network structures. For instance, these conditions can be realized in situations such as (I) when each particle is directly connected symmetrically to more than half of the population of particles, (II) when one or more leaders are directly connected to every particle, or (III) when the population is divided into three major groups with directed connections between them. Ultimately, these a priori conditions encompass network structures ranging from center-based connections to almost complete connections. Several concrete examples will be discussed in \Cref{sec_ex}.
\end{remark}

\subsection{Laplacian dynamics on strongly-connected temporal graphs}
In this subsection, we study the asymptotic consensus in Laplacian dynamics
\begin{align}\label{system3}
    \dot\sX(t)=-L(t)\sX(t),\quad\mbox{for}~~\sX\in\mathbb R^N,
\end{align}
where a temporal graph $\mathcal G(t)$, whose adjacency matrix is denoted by $(\gamma_{ij}(t))_{i,j}$, is strongly connected and satisfies the following a priori assumptions $(\mathcal B)$: for each $t\ge0$,

\vspace{.1cm}

\begin{enumerate}[label={($\mathcal B$\arabic*)}]
    \item there exist upper and lower bounds for number of neighbors, i.e.,
    \begin{align*}
        N_m\le\big|\{~j\in[N]\setminus\{i\}~:~\gamma_{ij}(t)\ne 0\}\big|\le N_M,\quad\forall~i\in[N],
    \end{align*}
    where $N_m\ge\lfloor N/2\rfloor$,
    \item there exist constants $0\le\epsilon$ and $0<\gamma_m$ such that $n\epsilon\le\gamma_m$,
    \begin{align}\label{assumB}
        \frac{\gamma_{ij}+\gamma_{ji}}{2}\in[\gamma_m,\infty)\quad\mbox{and}\quad\frac{\gamma_{ij}-\gamma_{ji}}{2}\in[-\epsilon,\epsilon],
    \end{align}
    where a positive constant $n>2$ satisfies 
    \begin{align}\label{n_cond}
        \frac{n-2}{s(n-1)}(N_m-s+1)-\frac{2}{(n+1)^2}\Big(3s^2-(N_M+2N_m+1)s+NN_M\Big)\ge\delta,
    \end{align}
    for all $1\le s\le\lfloor N/2\rfloor$ and some constant $\delta>0$,
    \item the edge set $\mathcal E(t)$ changes at most countable times.
\end{enumerate}
\begin{remark}\label{assumB_rmk}
    The relations \eqref{assumB} in $(\mathcal B2)$ imply that the graph must have symmetric edges but may exhibit possibly asymmetric weights. Moreover, as specified in $(\mathcal B1)$, a symmetric edge set $\mathcal E$, where each vertex has at least half of the other vertices as neighbors, guarantees that the graph is strongly-connected.
\end{remark}

\vspace{.1cm}

Note that the condition \eqref{assumB} implies that $\mathcal E(t)$ is symmetric and $\gamma_{ij}(t)\ge0$ for each $i,j\in[N]$ and $t\ge0$. Here, we introduce a fluctuation vector
\begin{align*}
    \hat\sX(t):=\sX(t)-\frac{1}{N}(\mathds{1}^{\mathsf{T}}X(t))\mathds{1}.
\end{align*}
Note that $\mathds{1}^{\mathsf{T}}\sX/N$ is the average value of particles $\{\sx_i\}_{i\in[N]}$ and it is not conserved in time as the graph is not symmetric.
\begin{theorem}\label{asymm_thm}
    The Laplacian dynamics \eqref{system3} with the assumptions $(\mathcal B)$ exhibits asymptotic consensus. Moreover, the distance to average exponentially decays,
    \begin{align*}
        \|\hat\sX(t)\|\le e^{-\frac{\gamma_m\delta t}{N}}\|\hat\sX(0)\|,\quad\forall~t\ge0.
    \end{align*}
\end{theorem}
\begin{proof}
    As in \Cref{a_priori_thm}, all estimates rely on the Carathéodory solution guaranteed by ($\mathcal B3$). The proof is based on a Lyapunov function
    \[G(t):=\frac{1}{2}\|\hat\sX(t)\|^2.\]
    For notation simplicity, we omit the time dependency. By definition of $\hat\sX$ and the linear ODE \eqref{system3}, one can get
    \begin{align*}
        \frac{d}{dt}G(t)&=\left(\sX-\frac{1}{N}(\mathds{1}^{\mathsf{T}}\sX)\mathds{1}\right)^{\mathsf{T}}\left(-L\sX+\frac{1}{N}(\mathds{1}^{\mathsf{T}}L\sX)\mathds{1}\right)\\
        &=-\sX^{\mathsf{T}}L\sX+\frac{1}{N}(\mathds{1}^{\mathsf{T}}L\sX)(\sX^{\mathsf{T}}\mathds{1})=-\sX^{\mathsf{T}}L\sX+\frac{1}{N}\mathds{1}^{\mathsf{T}}\sX\mathds{1}^{\mathsf{T}}L\sX\\
        &=-\hat\sX^{\mathsf{T}}L\sX=-\hat\sX^{\mathsf{T}}L\hat\sX-\hat\sX^{\mathsf{T}}L\left(\frac{1}{N}(\mathds{1}^{\mathsf{T}}\sX)\mathds{1}\right)=-\hat\sX^{\mathsf{T}}L\hat\sX.
    \end{align*}
    Rewrite it as element-wise,
    \begin{align*}
        \frac{d}{dt}G(t)&=-\sum_{i\in[N]}\hat\sx_i^2\sum_{j\ne i}\gamma_{ij}+\sum_{i\in[N]}\sum_{j\ne i}\hat\sx_i\hat\sx_j\gamma_{ij}=-\sum_{i\in[N]}\sum_{j\ne i}\gamma_{ij}\hat\sx_i(\hat\sx_i-\hat\sx_j)\\
        &=-\frac{1}{2}\sum_{i,j\in [N]}\left(\underbrace{\frac{\gamma_{ij}+\gamma_{ji}}{2}(\hat\sx_i-\hat\sx_j)^2+\frac{\gamma_{ij}-\gamma_{ji}}{2}(\hat\sx_i^2-\hat\sx_j^2)}_{=:\mathfrak C_{ij}(t)}\right).
    \end{align*}
    In other words,
    \begin{align}\label{dG}
        \frac{d}{dt}G(t)=-\frac{1}{2}\sum_{(i,j)\in\mathcal E(t)}\mathfrak C_{ij}(t).
    \end{align}
    Here, we distinguish edges $(i,j)$ and $(j,i)$. First, we discriminate when $\mathfrak C_{ij}$ is guaranteed to be nonnegative. If $(i,j)\notin\mathcal E(t)$, $\mathfrak C_{ij}=0$. For $(i,j)\in\mathcal E(t)$, we have
    \begin{align*}
        \mathfrak C_{ij}&=\frac{\gamma_{ij}+\gamma_{ji}}{2}(\hat\sx_i-\hat\sx_j)^2+\frac{\gamma_{ij}-\gamma_{ji}}{2}(\hat\sx_i^2-\hat\sx_j^2)\\
        &\ge \gamma_m(\hat\sx_i-\hat\sx_j)^2-\epsilon|\hat\sx_i^2-\hat\sx_j^2|\ge \gamma_m\left((\hat\sx_i-\hat\sx_j)^2-\frac{1}{n}|\hat\sx_i^2-\hat\sx_j^2|\right).
    \end{align*}
    By simple calculation, one can see that
    \begin{align*}
        \mathfrak C_{ij}\ge0~~\iff~~(\hat\sx_i,\hat\sx_j)\in\mathcal A_n^c,
    \end{align*}
    where the set $\mathcal A_n$ is defined by
    \begin{align*}
        \mathcal A_n:=\left\{(a,b)\in\mathbb R^2:ab>0~~\mbox{and}~~\frac{n-1}{n+1}<\frac{|a|}{|b|}<\frac{n+1}{n-1}\right\}.
    \end{align*}
    Below, we provide a rigorous lower bound for $\mathfrak C_{ij}$ for two cases.

    \vspace{.2cm}

    \noindent$\bullet$ {\bf(Non-positive product $\hat\sx_i\hat\sx_j\le0$)} Here, we consider when
    \begin{align*}
        (i,j)\in\mathcal E(t)\quad\mbox{and}\quad\hat\sx_i\hat\sx_j\le0.
    \end{align*}
    Without loss of generality, we assume $|\hat\sx_i|\ge|\hat\sx_j|\ge0$ to get
    \begin{align*}
        \mathfrak C_{ij}&\ge \gamma_m(\hat\sx_i-\hat\sx_j)^2-\epsilon(\hat\sx_i^2-\hat\sx_j^2)\ge \gamma_m\left((\hat\sx_i-\hat\sx_j)^2-\frac{1}{n}(\hat\sx_i^2-\hat\sx_j^2)\right)\\
        &\ge \gamma_m\left(\frac{n-1}{n}\left(\hat\sx_i-\frac{n\hat\sx_j}{n-1}\right)^2-\frac{1}{n(n-1)}|\hat\sx_j|^2\right)\ge \gamma_m\Big(\frac{n-1}{n}(|\hat\sx_i|^2+|\hat\sx_j|^2)-\frac{1}{n(n-1)}|\hat\sx_j|^2\Big).
    \end{align*}
    Therefore, we have
    \begin{align*}
        \mathfrak C_{ij}\ge \gamma_m\left(\frac{n-1}{n}(|\hat\sx_i|^2+|\hat\sx_j|^2)-\frac{1}{n(n-1)}\min\big(|\hat\sx_i|^2,|\hat\sx_j|^2\big)\right),
    \end{align*}
    more simply,
    \begin{align}\label{Anc_lb}
        \mathfrak C_{ij}\ge \frac{\gamma_m(n-2)}{n-1}(|\hat\sx_i|^2+|\hat\sx_j|^2).
    \end{align}

    \vspace{.2cm}

    \noindent$\bullet$ {\bf(Close ratio $(\hat\sx_i,\hat\sx_j)\in\mathcal A_n$)} Here, we consider the case $(\hat\sx_i,\hat\sx_j)\in\mathcal A_n$. The domain can be more specified. Every point $(\hat\sx_i,\hat\sx_j)\in\mathcal A_n$ satisfies
    \begin{align}\label{cond_A_n}
        \hat\sx_i\hat\sx_j>0,\quad\frac{|\hat\sx_i|}{|\hat\sx_j|}\in\left(\frac{n-1}{n+1},\frac{n+1}{n-1}\right),\quad|\hat\sx_i|,|\hat\sx_j|\le|\hat\sx_{\hat M}|,
    \end{align}
    where the notation $|\hat\sx_{\hat M}|$ is set to be
    \begin{align*}
        |\hat\sx_{\hat M}(t)|:=\max_{l\in[N]}|\hat\sx_l(t)|.
    \end{align*}
    Using \eqref{cond_A_n} and
    \begin{align*}
        |\hat\sx_i^2-\hat\sx_j^2|\le\frac{4n|\hat\sx_{\hat M}|^2}{(n+1)},
    \end{align*}
    one can derive
    \begin{align*}
        \mathfrak C_{ij}&\ge\frac{\gamma_{ij}-\gamma_{ji}}{2}(\hat\sx_i^2-\hat\sx_j^2)\ge-\epsilon|\hat\sx_i^2-\hat\sx_j^2|\ge-\frac{\gamma_m}{n}\min\left(\frac{4n|\hat\sx_{\hat M}|^2}{(n+1)^2},\max(\hat\sx_i^2,\hat\sx_j^2)\right).
    \end{align*}
    More simply, 
    \begin{align}\label{An_lb}
        \mathfrak C_{ij}\ge-\frac{4\gamma_m|\hat\sx_{\hat M}|^2}{(n+1)^2},\quad\forall~(\hat\sx_i,\hat\sx_j)\in\mathcal A_n.
    \end{align}

    \vspace{.2cm}

    \noindent Back to \eqref{dG}, decompose it as
    \begin{align*}
        \frac{d}{dt}G(t)\le-\frac{1}{2}\sum_{\substack{(i,j)\in\mathcal E(t)\\
        (\hat\sx_i,\hat\sx_j)\in\mathcal A_n}}\mathfrak C_{ij}(t)-\frac{1}{2}\sum_{\substack{(i,j)\in\mathcal E(t)\\
        (\hat\sx_i,\hat\sx_j)\in\mathcal A_n^c\\\hat\sx_i\hat\sx_j\le0}}\mathfrak C_{ij}(t).
    \end{align*}
    Now, we use the structure $(\mathcal B1)$ of the graph. Assume that $\|\hat\sX\|>0$. There there exists $s\le\lfloor N/2\rfloor$ such that
    \begin{align*}
        \begin{dcases}
            \hat\sx_i>0,\quad\forall~1\le i\le s,\\
            \hat\sx_i\le0,\quad\forall~s+1\le i\le N,
        \end{dcases}
    \end{align*}
    without loss of generality. Since $\mathds{1}^{\mathsf{T}}\hat\sX=0$, it is well defined. By the assumptions, each vertex $i$ has at least $N_m-(s-1)$ neighbors in $\{s+1,\cdots,N\}$. Thus, one has
    \begin{align}\label{dG_Anc}
        \begin{aligned}
            -\frac{1}{2}\sum_{\substack{(i,j)\in\mathcal E(t)\\
            (\hat\sx_i,\hat\sx_j)\in\mathcal A_n^c\\\hat\sx_i\hat\sx_j\le0}}\mathfrak C_{ij}(t)&\le-\frac{1}{2}\sum_{\substack{(i,j)\in\mathcal E(t)\\
            (\hat\sx_i,\hat\sx_j)\in\mathcal A_n^c\\\hat\sx_i\hat\sx_j\le0}}\frac{\gamma_m(n-2)}{n-1}(|\hat\sx_i|^2+|\hat\sx_j|^2)\\
            &\le-\frac{\gamma_m(n-2)}{n-1}\sum_{i=1}^s\sum_{\substack{j\ge s+1\\(i,j)\in\mathcal E(t)}}(|\hat\sx_i|^2+|\hat\sx_j|^2)\\
            &\le-\frac{\gamma_m(n-2)}{n-1}\big(N_m-(s-1)\big)\sum_{i=1}^s|\hat\sx_i|^2\\
            &\le-\frac{\gamma_m(n-2)}{s(n-1)}\Big(N_m-(s-1)\Big)|\hat\sx_{\hat M}|^2,
        \end{aligned}
    \end{align}
    where we use \eqref{Anc_lb} and 
    \begin{align*}
        |\hat\sx_{\hat M}|^2\le\left(\sum_{i=1}^s\hat\sx_i\right)^2\le s\sum_{i=1}^s|\hat\sx_i|^2.
    \end{align*}
    On the other hand, one can use \eqref{An_lb} to find
    \begin{align}\label{dG_An}
        \begin{aligned}
            -\frac{1}{2}\sum_{\substack{(i,j)\in\mathcal E(t)\\
            (\hat\sx_i,\hat\sx_j)\in\mathcal A_n}}\mathfrak C_{ij}(t)&\le\frac{1}{2}\sum_{\substack{(i,j)\in\mathcal E(t)\\
            (\hat\sx_i,\hat\sx_j)\in\mathcal A_n}}\frac{4\gamma_m|\hat\sx_{\hat M}|^2}{(n+1)^2}\\
            &\hspace{-2cm}\le\frac{2\gamma_m|\hat\sx_{\hat M}|^2}{(n+1)^2}\Bigg[s(s-1)+\Big((N-s)N_M-2s\big(N_m-(s-1)\big)\Big)\Bigg].
        \end{aligned}
    \end{align}
    Integrating \eqref{dG_Anc} and \eqref{dG_An} with \eqref{n_cond}, one can derive
    \begin{align*}
        \frac{d}{dt}G(t)&\le-\gamma_m|\hat\sx_{\hat M}|^2\left(\frac{n-2}{s(n-1)}(N_m-s+1)-\frac{2}{(n+1)^2}\Big(3s^2-(N_M+2N_m+1)s+NN_M\Big)\right)\\
        &\le-\frac{\gamma_m}{N}\|\hat\sX\|^2\left(\frac{n-2}{s(n-1)}(N_m-s+1)-\frac{2}{(n+1)^2}\Big(3s^2-(N_M+2N_m+1)s+NN_M\Big)\right)\\
        &\le-\frac{\gamma_m\delta}{N}\|\hat\sX\|^2\le0,
    \end{align*}
    where we use $\|\hat\sX\|\le N|\hat\sx_{\hat M}|^2$.
    Hence, we have
    \begin{align*}
        \frac{d}{dt}\|\hat\sX\|\le-\frac{\gamma_m\delta}{nN}\|\hat\sX\|,\quad\forall~t\ge0,
    \end{align*}
    which leads to the desired result.
\end{proof}
Now, we provide the sufficient framework for asymptotic flocking in \eqref{singular_adapCS} under suitable priori assumptions.
\begin{corollary}\label{cor_thm3}
    Suppose that constants $a_m\in(0,1)$, $\delta>0$ and two natural numbers $\lfloor N/2\rfloor\le N_m\le N_M<N$ such that
    \begin{align}\label{n2}
        n:=\frac{2a_m(N_m+1)}{N_M+1-a_m(N_m+1)}>2
    \end{align}
    and \eqref{n_cond} is satisfied. Assume that initial data $\{(\bx_i^0,\bv_i^0)\}_{i\in[N]}\in\mathcal D$ satisfy
    \begin{align*}
        0<a_m\le\min_{i,j\in[N]}\frac{\bv_i^0\cdot\bv_j^0}{|\bv_i^0||\bv_j^0|}
    \end{align*}
    and a temporal graph $\mathcal G(t)$ satisfies $(\mathcal B1)$. Then, the solution $\{(\bx_i,\bv_i)\}_{i\in[N]}$ exhibits asymptotic flocking.
\end{corollary}
\begin{proof}
    By \Cref{asymm_thm}, it is enough to show that there exists $n>2$ which satisfies the assumption $(\mathcal B2)$. Note that the system \eqref{singular_adapCS} is the Laplacian dynamics on the graph
    \begin{align*}
        \left(\frac{\kappa\psi_{ij}a_{ij}}{\mathcal I_i}\right)_{i,j}.
    \end{align*}
    The assumption $(\mathcal B1)$ and the definition of $\mathcal I_i$ and $\psi_{ij}$, one has
    \begin{align*}
        N_m+1\le \mathcal I_i(t)\le N_M+1,\quad\forall~i\in[N],~~\forall~t\ge0.
    \end{align*}
    In addition, from \Cref{a_bound}, we get
    \begin{align*}
        1\ge a_{ij}(t)\ge a_m,\quad\forall~(i,j)\in\mathcal E(t),~~\forall~t\ge0.
    \end{align*}
    Therefore, there is upper and lower bounds for weights in graph: if $\psi_{ij}=1$,
    \begin{align*}
        \gamma_m:=\frac{\kappa a_m}{N_M+1}\le\frac{\kappa a_{ij}}{\mathcal I_i}\le\frac{\kappa}{N_m+1}.
    \end{align*}
    Also, one can see the difference between $(i,j)$ and $(j,i)$ elements in the graph \eqref{CS_graph} as
    \begin{align*}
        \frac{\kappa}{2}\bigg|\frac{a_{ij}}{\mathcal I_i}-\frac{a_{ji}}{\mathcal I_j}\bigg|\le\frac{\kappa}{2}\frac{N_M+1-a_m(N_m+1)}{(N_M+1)(N_m+1)}=:\epsilon.
    \end{align*}
    The parameter values of $\gamma_m$ and $\epsilon$ satisfies $\gamma_m=n\epsilon$. Now, we again use the assumption to see that $(\mathcal B2)$ is satisfied, which leads to the end of the proof.
\end{proof}
\begin{remark}
    The a priori conditions in \Cref{cor_thm3} fundamentally depend on whether the parameters can produce sufficiently large values of $n$ on account of \eqref{n_cond}. This means that the setting \eqref{n2} requires that two parameters $N_m$ and $N_M$ must be sufficiently close to each other. Furthermore, it provides the following upper bound
    \begin{align*}
        n=\frac{2a_m(N_m+1)}{N_M+1-a_m(N_m+1)}\le\frac{2a_m}{1-a_m},
    \end{align*}
    which implies that the parameter $a_m$ need to be close to 1. In other words, these a priori conditions represent relatively strong assumptions. However, they correspond to scenarios that are indeed realizable within our model. The examples are depicted in \Cref{sec_ex}.
\end{remark}
\section{Examples}\label{sec_ex}

In this section, we provide several examples to gain a deeper understanding of the model and theorems presented in the previous section. First, we examine the sufficient initial conditions for a two-particle toy model to exhibit either flocking or dispersion. Although this model is not directly generalized to $N$-body systems in the paper, it implicitly suggests the asymptotic behavior of two pre-formed clusters. The second and third subsections provide concrete examples of graphs that satisfy the a priori assumptions in \Cref{a_priori_thm} and \Cref{asymm_thm}, respectively, illustrating how different network structures can fulfill the required conditions.

\subsection{Two-particles toy model}
Here, we explicitly present the sufficient conditions for asymptotic flocking and separation in the two-particles toy model, i.e., a system \eqref{singular_adapCS} with $N=2$. The motivation is to show the aggregation and separation of two groups, where the elements within each group interact in a manner similar to a single entity at the initial time.
\begin{lemma}\label{two_body}
    The solution $\{(\bx_i,\bv_i)\}_{i=1,2}$ to two-particles system \eqref{singular_adapCS} with initial data such that
    \begin{align*}
        \|\bx_1^0-\bx_2^0\|<d\quad\mbox{and}\quad a_m:=\frac{\bv_1^0\cdot\bv_2^0}{|\bv_1^0||\bv_2^0|}>0
    \end{align*}
    satisfies the following assertions:\newline
    \noindent(i) {\bf(Non-flocking)} If the initial data satisfy
    \begin{align*}
        \begin{dcases}
            \|\bx_1^0-\bx_2^0\|^2+\frac{1}{\kappa}(\bx_1^0-\bx_2^0)\cdot(\bv_1^0-\bv_2^0)+\frac{\|\bv_1^0-\bv_2^0\|^2}{2\kappa^2}>d^2,\quad&\mbox{if}~~(\bx_1^0-\bx_2^0)\cdot(\bv_1^0-\bv_2^0)\ge0,\\
            \|\bx_1^0-\bx_2^0\|^2+\frac{1}{\kappa a_m}(\bx_1^0-\bx_2^0)\cdot(\bv_1^0-\bv_2^0)+\frac{\|\bv_1^0-\bv_2^0\|^2}{2\kappa^2}>d^2,&\mbox{if}~~(\bx_1^0-\bx_2^0)\cdot(\bv_1^0-\bv_2^0)<0,
        \end{dcases}
    \end{align*}
    then the solution shows the separation of particles, i.e.,
    \[\lim_{t\to\infty}\|\bx_1(t)-\bx_2(t)\|=\infty,\quad\lim_{t\to\infty}\|\bv_1(t)-\bv_2(t)\|\ne0.\]
    \noindent(ii) {\bf(Flocking)} If the initial data satisfy
    \begin{align*}
        \begin{dcases}
            \|\bx_1^0-\bx_2^0\|^2+\frac{1}{\kappa a_m}(\bx_1^0-\bx_2^0)\cdot(\bv_1^0-\bv_2^0)+\frac{\|\bv_1^0-\bv_2^0\|^2}{2(a_m\kappa)^2}<d^2,\quad&\mbox{if}~~(\bx_1^0-\bx_2^0)\cdot(\bv_1^0-\bv_2^0)\ge0,\\
            \|\bx_1^0-\bx_2^0\|^2+\frac{1}{\kappa}(\bx_1^0-\bx_2^0)\cdot(\bv_1^0-\bv_2^0)+\frac{\|\bv_1^0-\bv_2^0\|^2}{2(a_m\kappa)^2}<d^2,&\mbox{if}~~(\bx_1^0-\bx_2^0)\cdot(\bv_1^0-\bv_2^0)<0,
        \end{dcases}
    \end{align*}
    where the constant $a_m$ is defined by
    \begin{align*}
        a_m:=\frac{\bv_1^0\cdot\bv_2^0}{|\bv_1^0||\bv_2^0|},
    \end{align*}
    then the solution attain the aggregation of particles, i.e.,
    \[\lim_{t\to\infty}\|\bx_1(t)-\bx_2(t)\|<d,\quad\lim_{t\to\infty}\|\bv_1(t)-\bv_2(t)\|=0.\]
\end{lemma}
\begin{proof}
    For notation simplicity, we denote
    \begin{align*}
        a(t):=\frac{\bv_1\cdot\bv_2}{|\bv_1||\bv_2|}\quad\mbox{and}\quad u(t):=(\bx_1-\bx_2)\cdot(\bv_1-\bv_2).
    \end{align*}
    From \eqref{singular_adapCS}, one has
    \begin{align}\label{spat_diff_deri}
        \frac{d}{dt}\|\bx_1-\bx_2\|^2=2u(t).
    \end{align}    
    We examine the derivative of spatial distance \eqref{spat_diff_deri}.
    \begin{align}\label{u_deri}
        \frac{d}{dt}u(t)&=-\kappa\psi_{12}a(t)u(t)+\|\bv_1-\bv_2\|^2.
    \end{align}
    The Hegselmann-Krause kernel $\psi_{12}=\chi_{[0,d)}(\|\bx_i-\bx_j\|)$ determines the connection between particles. Once they lose the connection, it can not be reconnected since $u(t)$ is nonnegative and monotone increasing from the lost time. Hence, for the constant $T\in(0,\infty]$ defined by
    \begin{align*}
        T:=\sup\{t>0~:~\|\bx_1(t)-\bx_2(t)\|<d\},
    \end{align*}
    we have
    \begin{align*}
        \|\bx_1(t)-\bx_2(t)\|<d,\quad\forall~0\le t<T.
    \end{align*}
    From now on, we restrict our computation to the interval $0\le t<T$, where $\psi_{12}$ is identically 1. Thus, it requires no further consideration and is omitted. Again from \eqref{singular_adapCS}, one has
    \begin{align*}
        \frac{d}{dt}\|\bv_1-\bv_2\|^2=-2\kappa a(t)\|\bv_1-\bv_2\|^2,
    \end{align*}
    which yields
    \begin{align}\label{v_diff}
        \|\bv_1(t)-\bv_2(t)\|^2=\|\bv_1^0-\bv_2^0\|^2\exp\left(-2\kappa\int_0^ta(s)ds\right).
    \end{align}
    Combining \eqref{u_deri} and \eqref{v_diff}, we get
    \begin{align*}
        u(t)=\exp\left(-\kappa\int_0^ta(s)ds\right)\left(u(0)+\|\bv_1^0-\bv_2^0\|^2\int_0^t\exp\left(-\kappa\int_0^sa(\tilde s)d\tilde s\right)ds\right).
    \end{align*}
    Using this solution $u(t)$, we can derive the explicit formula for spatial distance from \eqref{spat_diff_deri} 
    \begin{align}\label{spat_diff}
        \begin{aligned}
            &\|\bx_1(t)-\bx_2(t)\|^2-\|\bx_1^0-\bx_2^0\|^2\\
            &\hspace{.5cm}=\int_0^t\exp\left(-\kappa\int_0^{\tilde t}a(s)ds\right)\left(u(0)+\|\bv_1^0-\bv_2^0\|^2\int_0^{\tilde t}\exp\left(-\kappa\int_0^sa(\tilde s)d\tilde s\right)ds\right)d\tilde t.
        \end{aligned}
    \end{align}
    If there exists $t_*>0$ such that $\|\bx_1(t_*)-\bx_2(t_*)\|=d$, it is separation of particles ($T=t_*$), otherwise they are permanently connected over time ($T=\infty$). So, it is enough to see the asymptotic value of $\|\bx_1-\bx_2\|$ via \eqref{spat_diff}. Below, we focus on the case where the inner product of the relative distance and relative velocity at the initial time is nonnegative. The negative case is analogously and is therefore omitted.
    
    \noindent$(i)$ One can find the lower bound of asymptotic spatial distance: with $a(t)\le 1$,
    \begin{align*}
        \|\bx_1(t)-\bx_2(t)\|^2&\ge\|\bx_1^0-\bx_2^0\|^2+\int_0^te^{-\kappa\tilde t}\left(u(0)+\|\bv_1^0-\bv_2^0\|^2\int_0^{\tilde t}e^{-\kappa s}ds\right)dt\\
        &=\|\bx_1^0-\bx_2^0\|^2+\frac{u(0)}{\kappa}(1-e^{\kappa t})+\frac{\|\bv_1^0-\bv_2^0\|^2}{\kappa^2}\left(\frac{1}{2}-e^{-\kappa t}+\frac{e^{-2\kappa t}}{2}\right),
    \end{align*}
    hence,
    \begin{align*}
        \lim_{t\to\infty}\|\bx_1(t)-\bx_2(t)\|^2\ge\|\bx_1^0-\bx_2^0\|^2+\frac{u(0)}{\kappa}+\frac{\|\bv_1^0-\bv_2^0\|^2}{2\kappa^2}>d^2.
    \end{align*}
    In consequence, particles lose the connection at certain time $t_*\in(0,\infty)$. Since the velocity alignment cannot be reached in finite time by \eqref{v_diff}, i.e., $\bv_1(t_*)\ne\bv_2(t_*)$. For there is no acceleration after $t\ge t_*$, their velocities are fixed.
    
    \noindent$(ii)$
    Similarly, we use $a(t)\ge a(0)=a_m$ by \Cref{a_bound} to derive
    \begin{align*}
        \lim_{t\to\infty}\|\bx_1(t)-\bx_2(t)\|^2\le\|\bx_1^0-\bx_2^0\|^2+\frac{u(0)}{a_m\kappa}+\frac{\|\bv_1^0-\bv_2^0\|^2}{2(a_m\kappa)^2}\le d^2,
    \end{align*}
    with
    \begin{align*}
        \|\bx_1(t)-\bx_2(t)\|<d,\quad\forall~0<t<\infty.
    \end{align*}
    That is, their connection is permanent and by \eqref{v_diff}, one can obtain
    \begin{align*}
        \lim_{t\to\infty}\|\bv_1(t)-\bv_2(t)\|=0,
    \end{align*}
    which completes the proof.
\end{proof}
As mentioned in the beginning of this subsection, while the application of \Cref{two_body} remains unclear, it is expected to be useful to demonstrate the aggregation and separation of two groups.
\subsection{Examples of $(\mathcal A)$}\label{subsec_assumA}
In this subsection, we examine examples of graphs satisfying the assumption $(\mathcal A)$, which requires the time-varying graph to be strongly connected and neighbor-connected at each time. Specifically, we present examples for both undirected and directed graphs, focusing on their structure in terms of edges, while disregarding weights. Thus, the graph considered are of the form $\mathcal G=(\mathcal V,\mathcal E)$, with adjacency matrices denoted by $A\in\{0,1\}^{N\times N}$, where $N=|\mathcal V|$. If weights are included, additional asymmetry can be introduced to the graph.

\begin{figure}[h]
    \centering
    \begin{tikzpicture}[scale=0.5]
        \node at (-3.2,5) {$A_{\text{asym}}=$};
    
        \draw[thick] (-0.5,0) .. controls (-1.5,2.5) and (-1.5,7.5) .. (-0.5,10);
        \draw[thick] (10.5,0) .. controls (11.5,2.5) and (11.5,7.5) .. (10.5,10);
        
        \fill[blue!15] (0,5) -- (0,10) -- (5,10) -- (5,5) -- cycle;
        \fill[blue!15] (5,2) -- (5,5) -- (8,5) -- (8,2) -- cycle;
        \fill[blue!15] (8,0) -- (8,2) -- (10,2) -- (10,0) -- cycle;
        
        \draw[thick] (0,5) -- (0,10) -- (5,10) -- (5,5) -- cycle;
        \draw[thick] (5,2) -- (5,5) -- (8,5) -- (8,2) -- cycle;
        \draw[thick] (8,0) -- (8,2) -- (10,2) -- (10,0) -- cycle;
        \draw (0,0) -- (0,2) -- (5,2) -- (5,0) -- cycle;
        \draw (0,2) -- (0,5) -- (5,5) -- (5,2) -- cycle;
        \draw (5,0) -- (5,2) -- (8,2) -- (8,0) -- cycle;
        \draw (5,5) -- (5,10) -- (8,10) -- (8,5) -- cycle;
        \draw (8,5) -- (8,10) -- (10,10) -- (10,5) -- cycle;
        \draw (8,2) -- (8,5) -- (10,5) -- (10,2) -- cycle;
    
        \draw[<->] (0,10.5) -- (5,10.5);
        \draw[<->] (5,-0.5) -- (8,-0.5);
        \draw[<->] (8,10.5) -- (10,10.5);
        
        \node at (2.5,7.5) {1};
        \node at (6.5,7.5) {1};
        \node at (9,7.5) {0};
        \node at (2.5,3.5) {0};
        \node at (6.5,3.5) {1};
        \node at (9,3.5) {1};
        \node at (2.5,1) {1};
        \node at (6.5,1) {0};
        \node at (9,1) {1};
        \node at (2.5,11) {$\mathcal V_1$};
        \node at (6.5,-1) {$\mathcal V_2$};
        \node at (9,11) {$\mathcal V_3$};
    \end{tikzpicture}
    \caption{Adjacency matrix representation of a {\it directed} neighbor-connected graph, as described in \Cref{subsec_assumA}. The graph illustrates the neighbor-connectivity defined in \Cref{def_neighbor_conn}.}
    \label{ex_asym_nbcnt}
\end{figure}
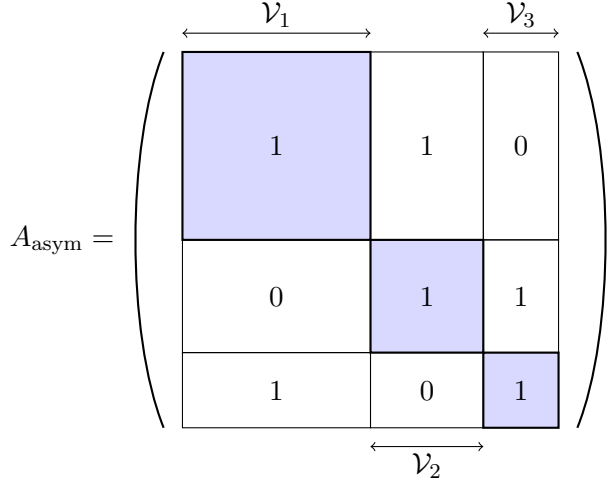

For the undirected graph case, we consider the following examples:
\begin{itemize}
	\item A singer-leader graph, denoted by $\mathcal G_{N,1}^{(\text{nc})}$, where a designated vertex $v_{\text{lead}}$ is connected to all other Vertices, and no additional edges are present.
	\item A two-leaders graph, denoted by $\mathcal G_{N,2}^{(\text{nc})}$, where two Vertices $v_{\text{lead},1}$ and $v_{\text{lead},2}$, each connected to all other vertices, with no other edges.
\end{itemize}
\begin{align}\label{ncG_adj}
    A^{\text{(nc)}}_{N,1}:=
    \begin{bmatrix}
        0 & 1 & 1 & \cdots & 1 \\
        1 & 0 & 0 & \cdots & 0 \\
        1 & 0 & 0 & \cdots & 0 \\
        \vdots & \vdots & \vdots & \ddots & \vdots \\
        1 & 0 & 0 & \cdots & 0
    \end{bmatrix},\quad
    A^{\text{(nc)}}_{N,2}:=
    \begin{bmatrix}
        0 & 1 & \cdots & 1 & 0 \\
        1 & 0 & \cdots & 0 & 1 \\
        \vdots & \vdots & \vdots & \ddots & \vdots \\
        1 & 0 & \cdots & 0 & 1 \\
        0 & 1 & \cdots & 1 & 0
    \end{bmatrix}.
\end{align}
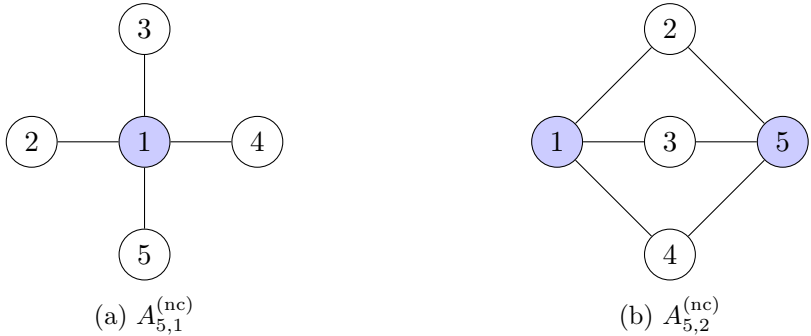
\begin{figure}[h]
    \centering
    \begin{subfigure}[t]{0.45\textwidth}
        \centering
        \begin{tikzpicture}[scale=1, every node/.style={circle, draw, fill=none}]
            \node (1)[draw, fill=blue!20] at (0, 0) {1};
            \node (2) at (-1.5, 0) {2};
            \node (3) at (0, 1.5) {3};
            \node (4) at (1.5, 0) {4};
            \node (5) at (0, -1.5) {5};
            
            \draw (1) -- (2);
            \draw (1) -- (3);
            \draw (1) -- (4);
            \draw (1) -- (5);
        \end{tikzpicture}
        \caption{$A^{\text{(nc)}}_{5,1}$}
    \end{subfigure}
    \begin{subfigure}[t]{0.45\textwidth}
        \centering
        \begin{tikzpicture}[scale=1, every node/.style={circle, draw, fill=none}]
            \node (1)[draw, fill=blue!20] at (-1.5, 0) {1};
            \node (2) at (0, 1.5) {2};
            \node (3) at (0, 0) {3};
            \node (4) at (0, -1.5) {4};
            \node (5)[draw, fill=blue!20] at (1.5, 0) {5};
            
            \draw (1) -- (2);
            \draw (1) -- (3);
            \draw (1) -- (4);
            \draw (2) -- (5);
            \draw (3) -- (5);
            \draw (4) -- (5);
        \end{tikzpicture}
        \caption{$A^{\text{(nc)}}_{5,2}$}
    \end{subfigure}
    \caption{Examples of neighbor-connected graphs}
    \label{graph_vertex5}
\end{figure}
\noindent The adjacency matrices for these examples are provided in \eqref{ncG_adj}. The graphical representations of these graphs with a vertex set of size 5 are shown in \Cref{graph_vertex5}. For the case of directed graphs, we construct an example by partitioning the vertex set into three disjoint, nonempty groups, denoted by $\mathcal V_1$, $\mathcal V_2$ and $\mathcal V_3$. Within each subset, the induced subgraph is a directed complete graph, ensuring that the shortest path between any two vertices in the subset has a length 1. Inter-group edges are added directly as follows: $\mathcal V_1\to\mathcal V_2\to\mathcal V_3\to\mathcal V_1$. This construction satisfies the given conditions, with the adjacency matrix illustrated in \Cref{ex_asym_nbcnt}.
            
            
\subsection{Examples of $(\mathcal B)$}\label{subsec_assumB}
In this subsection, we provide examples of graphs satisfying the assumption $(\mathcal B)$. As we discussed in \Cref{assumB_rmk}, we consider a graph $\mathcal G=(\mathcal V,\mathcal E,\mathcal W)$ where the edge set $\mathcal E$ is symmetric. The construction process involves satisfying $(\mathcal B1)$ for the edge set and $(\mathcal B2)$ for the weights. To begin, we consider the edge set under the constraints of the first assumption. This can be categorized into two cases:
\begin{itemize}
    \item When $N_m=N_M\equiv n$, meaning that each vertex has exactly $n$ neighbors, a vertex-transitive $n$-regular graph is specified.
    \item In the case of $N_m\ne N_M$, i.e., the number of neighbors per vertex varies within a nontrivial range, a graph can be constructed if the following conditions are satisfied:
    \begin{align}\label{NmNM_rel}
        \frac{2}{3}(N-1)\ge N_m,\quad 2(N-1-N_m)\ge N_M
    \end{align}
    The construction proceeds as follows:
    \begin{enumerate}[label=(\roman*),itemsep=0em]
        \item Select $N_M+1$ vertices to form a complete subgraph.
        \item Independently, select $N_m+1$ vertices, ensuring that all vertices not included in the first subgraph are selected for the second subgraph, which also forms a complete subgraph. Any edges appear in both subgraphs are removed.
    \end{enumerate}
\end{itemize}
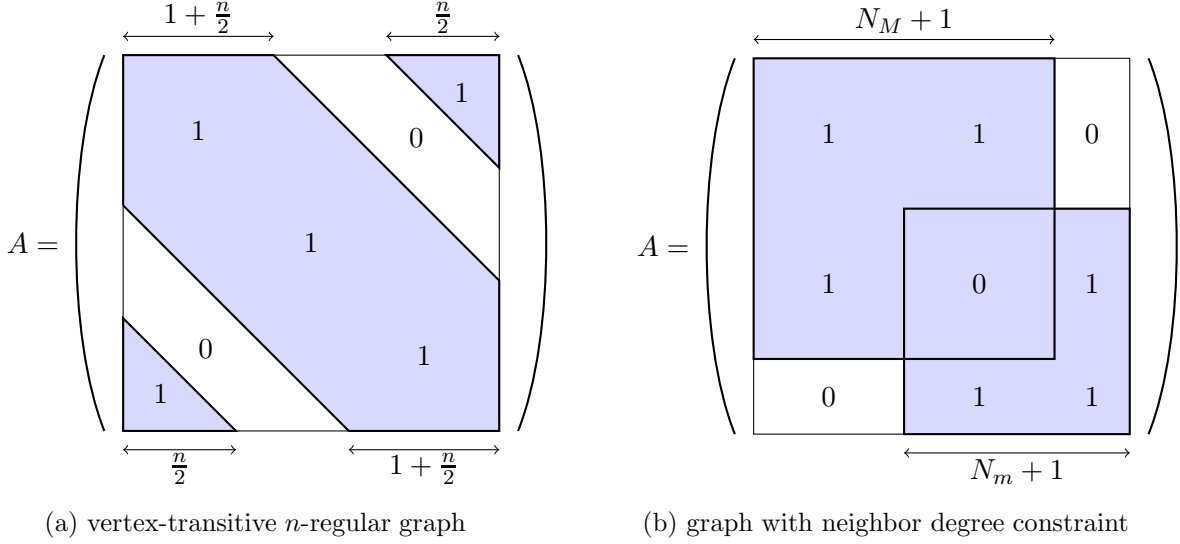
\begin{figure}
    \centering
    \begin{subfigure}{0.45\textwidth}
        \centering
        \begin{tikzpicture}[scale=0.5]
            \node at (-2.4,5) {$A=$};

            \draw[thick] (-0.5,0) .. controls (-1.5,2.5) and (-1.5,7.5) .. (-0.5,10);
            \draw[thick] (10.5,0) .. controls (11.5,2.5) and (11.5,7.5) .. (10.5,10);
        
            \fill[blue!15] (0,6) -- (0,10) -- (4,10) -- (10,4) -- (10,0) -- (6,0) -- cycle;
            \fill[blue!15] (7,10) -- (10,10) -- (10,7) -- cycle;
            \fill[blue!15] (0,0) -- (0,3) -- (3,0) -- cycle;
            
            \draw[thick] (0,6) -- (0,10) -- (4,10) -- (10,4) -- (10,0) -- (6,0) -- (0,6);
            \draw[thick] (7,10) -- (10,10) -- (10,7) -- cycle;
            \draw[thick] (0,0) -- (0,3) -- (3,0) -- cycle;
            \draw (0,0) -- (0,10) -- (10,10) -- (10,0) -- cycle;
        
            \draw[<->] (0,10.5) -- (4,10.5);
            \draw[<->] (7,10.5) -- (10,10.5);
            \draw[<->] (6,-0.5) -- (10,-0.5);
            \draw[<->] (0,-0.5) -- (3,-0.5);
            
            \node at (5,5) {1};
            \node at (8,2) {1};
            \node at (2,8) {1};
            \node at (7.8,7.8) {0};
            \node at (2.2,2.2) {0};
            \node at (1,1) {1};
            \node at (9,9) {1};
            \node at (2,11) {$1+\frac{n}{2}$};
            \node at (8,-1) {$1+\frac{n}{2}$};
            \node at (8.5,11) {$\frac{n}{2}$};
            \node at (1.5,-1) {$\frac{n}{2}$};
        \end{tikzpicture}
        \caption{vertex-transitive $n$-regular graph}
    \end{subfigure}
    \hfill
    \begin{subfigure}{0.45\textwidth}
        \centering
        \begin{tikzpicture}[scale=0.5]
            \node at (-2.4,5) {$A=$};
        
            \draw[thick] (-0.5,0) .. controls (-1.5,2.5) and (-1.5,7.5) .. (-0.5,10);
            \draw[thick] (10.5,0) .. controls (11.5,2.5) and (11.5,7.5) .. (10.5,10);
        
            \fill[blue!15] (0,2) -- (8,2) -- (8,10) -- (0,10) -- cycle;
            \fill[blue!15] (4,0) -- (10,0) -- (10,6) -- (4,6) -- cycle;
            
            \draw[thick] (0,2) -- (8,2) -- (8,10) -- (0,10) -- cycle;
            \draw[thick] (4,0) -- (10,0) -- (10,6) -- (4,6) -- cycle;
            \draw (0,0) -- (0,2) -- (4,2) -- (4,0) -- cycle;
            \draw (8,6) -- (8,10) -- (10,10) -- (10,6) -- cycle;
        
            \draw[<->] (0,10.5) -- (8,10.5);
            \draw[<->] (4,-0.5) -- (10,-0.5);
            
            \node at (2,1) {0};
            \node at (2,4) {1};
            \node at (2,8) {1};
            \node at (6,1) {1};
            \node at (6,4) {0};
            \node at (6,8) {1};
            \node at (9,1) {1};
            \node at (9,4) {1};
            \node at (9,8) {0};
            \node at (4,11) {$N_M+1$};
            \node at (7,-1) {$N_m+1$};
        
        \end{tikzpicture}
        \caption{graph with neighbor degree constraint}
    \end{subfigure}
    \caption{Adjacency matrix of the edge set $\mathcal E$ for a connected graph $([N], \mathcal E, \mathcal W)$ satisfying the neighbor degree constraint ($\mathcal B1$), with the weight information $\mathcal W$ excluded. The graph is connected, and each node connects to a limited number of neighbors.}
    \label{assumB_adj}
\end{figure}
These cases and their corresponding adjacency matrices are illustrated in \Cref{assumB_adj}. Next, we turn to the weights under the assumption $(\mathcal B2)$. Essentially, it is sufficient to introduce small fluctuations around $\gamma_m$. For example, one could set the weights as
\[(w_{ij},w_{ji})=(\gamma_m-\delta_1,\gamma_m+\delta_2),\]
with $\delta_1,\delta_2$ being small random values in $[-\epsilon,\epsilon]$ satisfying $\delta_1\delta_2\ge0$.
\section{Numerical results}\label{sec_num}
In this section, we present numerical simulations of Laplacian dynamics on a various types of graphs and adaptive CS dynamics. For all simulations, we use a 4th-order Runge-Kutta scheme with a time step $\Delta t=10^{-3}$.
\subsection{Laplacian dynamics on graphs}
This subsection provides numerical simulations of the Laplacian dynamics
\begin{align*}
    \dot\sX=-L_{\mathcal G}(t)\sX,
\end{align*}
on time-varying graphs $\mathcal G(t)$ discussed in Subsection \ref{subsec_assumA} and \ref{subsec_assumB}. The simulations focus on the evolution of the diameter $D(t)$ of the solution $\sX(t)=[\sx_1,\cdots,\sx_N]^{\mathsf{T}}$, defined as
\begin{align*}
    D(t):=\max_{i,j\in[N]}|\sx_i(t)-\sx_j(t)|.
\end{align*}
The simulations are conducted in three parts: (i) graphs satisfying assumption $(\mathcal A)$ and undirected (see \eqref{ncG_adj}), (ii) graphs satisfying assumption $(\mathcal A)$ and directed (see \Cref{ex_asym_nbcnt}), and (iii) graphs satisfying assumption B (see \Cref{assumB_adj}).
\begin{figure}
    \centering
    \begin{subfigure}{0.32\textwidth}
        \centering
        \begin{overpic}
            [width=1.1\linewidth]{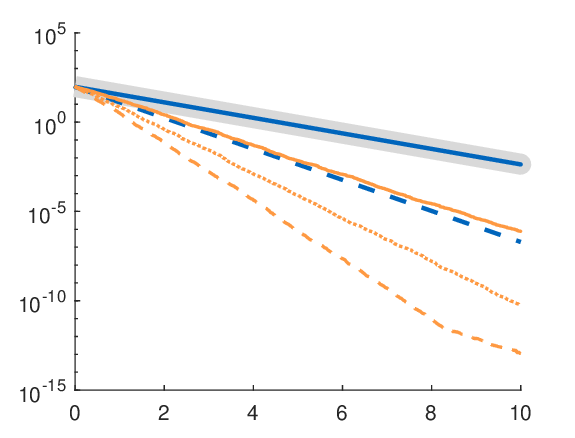}
            \put(-4,19){\rotatebox{90}{\small diameter $D(t)$}}
            \put(43,-2.3){\small time $t$}
        \end{overpic}
        \caption{undirected graph with $(\mathcal A)$}
        \label{sim_thm2_1}
    \end{subfigure}
    \begin{subfigure}{0.32\textwidth}
        \centering
        \begin{overpic}
            [width=1.1\linewidth]{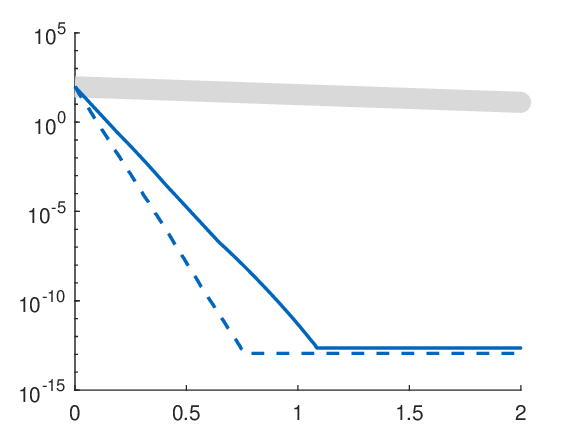}
        \end{overpic}
        \caption{directed graph with $(\mathcal A)$}
        \label{sim_thm2_2}
    \end{subfigure}
    \begin{subfigure}{0.32\textwidth}
        \centering
        \begin{overpic}
            [width=1.1\linewidth]{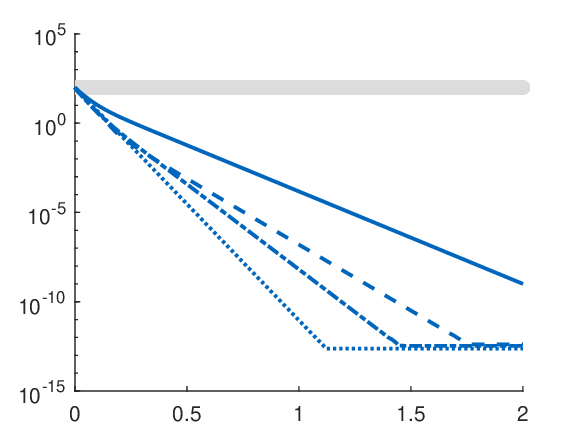}
        \end{overpic}
        \caption{graph with $(\mathcal B)$}
        \label{sim_thm3}
    \end{subfigure}
    \caption{Diameter evolution over time for various graph types. (a) Fixed one-leader (blue solid), fixed two-leaders (blue dashed), temporal one-leader (orange solid), temporal two-leaders (orange dashed), and temporal one or two leaders (orange dotted). (b) Fixed three groups graph (solid) and temporal three groups graph (dashed). (c) Graphs with different configurations: $(N_m, N_M) = (30, 30)$ (solid), $(N_m, N_M) = (30, 45)$ (dashed), $(N_m, N_M) = (30, 58)$ (dotted), and $(N_m, N_M) = (38, 38)$ (dash-dot). All experiments use a time-step $0.001$ with the same initial data generated by $\bx_i(0)\sim\mathcal U[0,100]$ for all experiments. The thick gray lines represent the upper bound estimated in \Cref{a_priori_thm} and \Cref{asymm_thm}, which are consistently larger than all experimental results, confirming the validity of our approach.}
    \label{sim_thm23}
\end{figure}

The simulations involve graphs with 60 vertices ($|\mathcal V|=60$), and results are illustrated in \Cref{sim_thm23}. The thick gray line represents the upper bounds provided in \Cref{a_priori_thm} and \Cref{asymm_thm}. While the decay constants are either very similar or extremely small, making the lines appear as a single one, they are, in fact, distinct lines derived from the settings of each graph. The graphs are plotted on a logarithmic $y$-axis, and all curves appear linear in their asymptotic regime. This indicates that their convergence rates decrease exponentially, differing only in the rate constants. This observation aligns with the results of our main results (\Cref{a_priori_thm}, \Cref{asymm_thm}), which also predict exponential decay.

\Cref{sim_thm2_1} illustrates the diameter evolution of solutions by two fixed graphs (blue) and three time-varying graphs (orange). The fixed cases include a one-leader graph and a two-leaders graph, while the time-varying cases consist of switching one-leader, switching two-leaders, and a mixed case alternating between the two in a half-and-half ratio. The results reveal a tendency for faster convergence when leaders are alternating within the same structure. This occurs because leader switching allows the temporal graph to provide more direct communication opportunities between vertices over time. As expected, the one-leader graph exhibits slower convergence than the two-leaders graph, as it has fewer edges $(N-1)$ compared to the two-leaders case $2(N-2)$. Furthermore, the one-leader graph reduces the estimate in \eqref{dich_estimate} to equalities except for $(N-1)$ rows relating to leader vertex, signifying it as the case closest to the slowest scenario among all graphs satisfying the assumption $(\mathcal A)$.

For the directed case, \Cref{sim_thm2_2} presents results for both fixed and time-varying graphs, with the latter constructed by randomly partitioning the vertex set $\mathcal V$ at each step. As in the undirected case, randomly changing structures guarantee faster convergence rates. However, even in the fixed case, there are instances where \eqref{dich_estimate} provides a relatively large value as a lower bound, leading to faster convergence rates than the upper bound suggested in \Cref{a_priori_thm}.
\begin{table}[h!]
\centering
\begin{tabular}{c||c|c}
\noalign{\hrule height 0.4mm}
$(N_m,N_M)$ & $\delta$ & $n$ \\
\noalign{\hrule height 0.4mm}
(30,30) & 1.2469e-04 & 326 \\
\hline
(30,45) & 9.6546e-5 & 365 \\
\hline
(30,58) & 1.2901e-4 & 396 \\
\hline
(38,38) & 0.0029 & 101 \\
\hline
\end{tabular}
\caption{Values of $\delta$ and $n$ used in \Cref{sim_thm3}}
\label{dlta_n}
\end{table}
Finally, \Cref{sim_thm3} examines scenarios with various $(N_m,N_M)$ satisfying \eqref{NmNM_rel}, consisting of two $n$-regular graph and two non-trivial graphs described in \Cref{assumB_adj}. For each case, $\delta$ and $n$ in $(\mathcal B2)$ are chosen as in \Cref{dlta_n}. The three cases with $N_m=30$ exhibits faster decay as $N_M$ increases. Even when $N_M$ is as large as 58, which means that some vertices directly connect to all but one vertex, the case of the 38-regular graph shows faster convergence. This suggests that the minimum number of neighbors significantly influences the convergence rate.

\subsection{Dynamics in the adaptive Cucker-Smale model}

Another important strength of the adaptive CS model lies in its ability to form {\it meaningful} clusters representing distinctly different opinions. Unlike clustering which emerges merely from disconnections caused by spatial distance, the clustering observed in the adaptive CS dynamics results from a clear divergence into several distinct orientations, which subsequently leads to network disconnections. To ensure a fair comparison, we configure both the adaptive CS model and the classical CS model with a connection kernel $\psi$ and normalizers $\mathcal I_i$ defined in \eqref{I_Psi} and examine their behavior under randomly selected initial conditions as follows: with $N=60$,
\begin{align*}
    \bx_i^0=0,\quad\bv_i^0\sim\mathcal U([-0.01,0.01]^2),\quad\kappa_{ij}^0=\kappa_{ji}^0\sim\mathcal U([-0.5,0.5]),\quad\forall~i,j\in[N].
\end{align*}
\begin{table}[h!]
    \centering
    \begin{tabular}{c|c|c|c|c}
        \noalign{\hrule height 0.4mm}
        \cellcolor{yellow!50}Case 1 & Case 2 & Case 3 & Case 4 & Case 5\\
        \noalign{\hrule height 0.4mm}
        \makecell{23 agents, $16.4^{\circ}$\\37 agents, $182.9^{\circ}$\\
        \begin{tikzpicture}
            \draw[thick] (0,0) circle(0.8);
            \draw[->,thick] (0,0) -- (0.9593,0.2823);
            \draw[->,thick] (0,0) -- (-0.9987,-0.0506);
        \end{tikzpicture}
        }&\makecell{12 agents, $139.4^{\circ}$\\48 agents, $306.3^{\circ}$\\
        \begin{tikzpicture}
            \draw[thick] (0,0) circle(0.8);
            \draw[->,thick] (0,0) -- (-0.7593,0.6508);
            \draw[->,thick] (0,0) -- (0.5920,-0.8059);
        \end{tikzpicture}
        }&\makecell{28 agents, $39.4^{\circ}$\\32 agents, $218.5^{\circ}$\\
        \begin{tikzpicture}
            \draw[thick] (0,0) circle(0.8);
            \draw[->,thick] (0,0) -- (0.7727,0.6347);
            \draw[->,thick] (0,0) -- (-0.7826,-0.6225);
        \end{tikzpicture}
        }&\makecell{29 agents, $100.7^{\circ}$\\3 agents, $157.3^{\circ}$\\28 agents, $291.3^\circ$\\
        \begin{tikzpicture}
            \draw[thick] (0,0) circle(0.8);
            \draw[->,thick] (0,0) -- (-0.1857,0.9826);
            \draw[->,thick] (0,0) -- (-0.9225,0.3859);
            \draw[->,thick] (0,0) -- (0.3633,-0.9317);
        \end{tikzpicture}
        }&\makecell{29 agents, $42.9^{\circ}$\\6 agents, $128.1^{\circ}$\\25 agents, $238.3^{\circ}$\\
        \begin{tikzpicture}
            \draw[thick] (0,0) circle(0.8);
            \draw[->,thick] (0,0) -- (0.7325,0.6807);
            \draw[->,thick] (0,0) -- (-0.6170,0.7869);
            \draw[->,thick] (0,0) -- (-0.5255,-0.8508);
        \end{tikzpicture}
        }\\
        \noalign{\hrule height 0.4mm}
        \cellcolor{yellow!50}Case 6 & Case 7 & Case 8 & Case 9 & \cellcolor{yellow!50}Case 10\\
        \noalign{\hrule height 0.4mm}
        \cellcolor{yellow!50}\makecell{26 agents, $4.7^{\circ}$\\31 agents, $170.7^{\circ}$\\3 agents, $286.4^{\circ}$\\
        \begin{tikzpicture}
            \draw[thick] (0,0) circle(0.8);
            \draw[->,thick] (0,0) -- (0.9932,0.1167);
            \draw[->,thick] (0,0) -- (-0.9869,0.1616);
            \draw[->,thick] (0,0) -- (0.2823,-0.9593);
        \end{tikzpicture}
        }&\makecell{32 agents, $157.6^{\circ}$\\28 agents, $346.3^{\circ}$\\
        \begin{tikzpicture}
            \draw[thick] (0,0) circle(0.8);
            \draw[->,thick] (0,0) -- (-0.9245,0.3811);
            \draw[->,thick] (0,0) -- (0.9715,-0.2368);
        \end{tikzpicture}
        }&\makecell{22 agents, $36.3^{\circ}$\\11 agents, $49.8^{\circ}$\\7 agents, $196.9^{\circ}$\\20 agents, $221.5^{\circ}$\\
        \begin{tikzpicture}
            \draw[thick] (0,0) circle(0.8);
            \draw[->,thick] (0,0) -- (0.4446,0.8957);
            \draw[->,thick] (0,0) -- (0.6455,0.7638);
            \draw[->,thick] (0,0) -- (-0.9568,-0.2907);
            \draw[->,thick] (0,0) -- (-0.7490,-0.6626);
        \end{tikzpicture}
        }&\makecell{29 agents, $82.1^{\circ}$\\31 agents, $260.3^{\circ}$\\
        \begin{tikzpicture}
            \draw[thick] (0,0) circle(0.8);
            \draw[->,thick] (0,0) -- (0.1374,0.9905);
            \draw[->,thick] (0,0) -- (-0.1685,-0.9857);
        \end{tikzpicture}
        }&\makecell{6 agents, $9.4^{\circ}$\\ 21 agents, $45.3^{\circ}$\\16 agents, $169.1^{\circ}$,\\16 agents, $223.6^{\circ}$\\
        \begin{tikzpicture}
            \draw[thick] (0,0) circle(0.8);
            \draw[->,thick] (0,0) -- (0.9866,0.1633);
            \draw[->,thick] (0,0) -- (0.7034,0.7108);
            \draw[->,thick] (0,0) -- (-0.9820,0.1891);
            \draw[->,thick] (0,0) -- (-0.7242,-0.6896);
        \end{tikzpicture}
        }\\
        \hline
    \end{tabular}
    \caption{(Adaptive CS model) Number of agents and the polar angle of the velocity vector for the largest clusters in the population at terminal time $t = 300$. Each cell shows the particle count and the corresponding polar angle for the largest clusters in that experiment. Below, the individual cluster angles are visualized with arrows, providing a geometric representation of the cluster orientations. Cases 1, 6, and 10 are marked to indicate that they are compared with the original CS model in \Cref{tab_CS}. Case 6 is additionally marked to indicate that a time snapshot is shown in \Cref{fig_snapshot_rng6}}
    \label{tab_adapCS}
\end{table}

We present the computational results for asymptotic group formation and their opinions by the adaptive CS dynamics with radius $d=1$ starting from 10 randomly assigned initial conditions in \Cref{tab_adapCS}. Here, the number of agents in the largest groups and their velocity vectors’ polar angles are recorded at $t=300$, which can be considered sufficient for group formation and intra-group velocity synchronization. The recorded major groups, each encompassing a substantial portion of the population, are meaningful clusters due to their size. As anticipated, the velocity vectors of these groups exhibit sufficient angular separation to indicate clustering.
\begin{table}[h!]
    \centering
    \begin{tabular}{c||c|c|c}
        \noalign{\hrule height 0.4mm}
        & $d=0.009$ & $d=0.01$ & $d=0.012$ \\
        \noalign{\hrule height 0.4mm}
        \cellcolor{yellow!50}Case 1
        &\makecell{24 agents, $57.2^{\circ}$\\21 agents, $150.1^{\circ}$\\14 agents, $161.7^{\circ}$ \\
        \begin{tikzpicture}
            \draw[thick] (0,0) circle(0.8);
            \draw[->,thick] (0,0) -- (0.5417,0.846);
            \draw[->,thick] (0,0) -- (-0.8669,0.4985);
            \draw[->,thick] (0,0) -- (-0.9494,0.3140);
        \end{tikzpicture}}&\makecell{24 agents, $86.7^{\circ}$\\35 agents, $148.9^{\circ}$\\
        \begin{tikzpicture}
            \draw[thick] (0,0) circle(0.8);
            \draw[->,thick] (0,0) -- (0.0576,0.9983);
            \draw[->,thick] (0,0) -- (-0.8563,0.5165);
        \end{tikzpicture}
        }&\makecell{60 agents, $134.4^{\circ}$\\
        \begin{tikzpicture}
            \draw[thick] (0,0) circle(0.8);
            \draw[->,thick] (0,0) -- (-0.6997,0.7145);
        \end{tikzpicture}
        }\\
        \hline
        \cellcolor{yellow!50}Case 6&\cellcolor{yellow!50}\makecell{11 agents, $61.0^{\circ}$\\23 agents, $67.8^{\circ}$\\25 agents, $120.6^{\circ}$\\
        \begin{tikzpicture}
            \draw[thick] (0,0) circle(0.8);
            \draw[->,thick] (0,0) -- (0.4848,0.8746);
            \draw[->,thick] (0,0) -- (0.3778,0.9259);
            \draw[->,thick] (0,0) -- (-0.5090,0.8607);
        \end{tikzpicture}
        }&\makecell{34 agents, $89.9^{\circ}$\\26 agents, $100.7^{\circ}$\\
        \begin{tikzpicture}
            \draw[thick] (0,0) circle(0.8);
            \draw[->,thick] (0,0) -- (0.0017,1);
            \draw[->,thick] (0,0) -- (-0.1857,0.9826);
        \end{tikzpicture}
        }&\makecell{60 agents, $92.5^{\circ}$\\
        \begin{tikzpicture}
            \draw[thick] (0,0) circle(0.8);
            \draw[->,thick] (0,0) -- (-0.0436,0.9990);
        \end{tikzpicture}
        }\\
        \hline
        \cellcolor{yellow!50}Case 10
        &\makecell{25 agents, $196.0^{\circ}$\\15 agents, $208.9^{\circ}$\\20 agents, $215.1^{\circ}$\\
        \begin{tikzpicture}
            \draw[thick] (0,0) circle(0.8);
            \draw[->,thick] (0,0) -- (-0.9613,-0.2756);
            \draw[->,thick] (0,0) -- (-0.8755,-0.4833);
            \draw[->,thick] (0,0) -- (-0.8181,-0.5750);
        \end{tikzpicture}
        }&\makecell{60 agents, $260.4^{\circ}$\\
        \begin{tikzpicture}
            \draw[thick] (0,0) circle(0.8);
            \draw[->,thick] (0,0) -- (-0.1668,-0.9860);
        \end{tikzpicture}
        }&\makecell{60 agents, $207.6^{\circ}$\\
        \begin{tikzpicture}
            \draw[thick] (0,0) circle(0.8);
            \draw[->,thick] (0,0) -- (-0.8862,-0.4366);
        \end{tikzpicture}
        }\\
        \hline
    \end{tabular}
    \caption{(Original CS model) Number of agents and the polar angle of the velocity vector for the largest clusters in the population at terminal time $t = 300$. The initial data are the same as in Cases 1, 6, and 10 from the experiments for the adaptive CS model (\Cref{tab_adapCS}). In this experiment, the asymptotic behavior of the system is investigated for three values of the parameter $d=0.009,0.01,0.012$. This setup allows for a comparison of the asymptotic behavior between the Original and Adaptive CS models. Case 6 with $d=0.009$ is marked to indicate that a time snapshot is shown in \Cref{fig_snapshot_rng6_orig}.}
    \label{tab_CS}
\end{table}

In contrast, the classical CS model with connection kernel and normalizers,
\begin{align}\label{CS_d}
    \dot\bx_i=\bv_i,\quad\dot\bv_i=\frac{1}{\mathcal I_i}\sum_{j=1}^N\frac{\psi_{ij}}{(1+\|\bx_i-\bx_j\|^2)^{1/2}}(\bv_j-\bv_i),
\end{align}
exhibits different clustering behavior. Due to its tendency for all velocity vectors to synchronize, the clustering process in the classical CS model is primarily driven by disconnections caused by spatial distance. In \Cref{tab_CS}, we present numerical results for the clustering and velocity vectors' polar angles for the classical CS model \eqref{CS_d} under the same initial conditions as Case 1, 6, and 10 in \Cref{tab_adapCS}, with varying radius $d$. Depending on the radius $d$, the system tends to either synchronize in a single direction or split into multiple groups with relatively small differences in velocity vector orientations. A clear observation is that the classical CS model results in either a large mono-cluster, a dispersed state of very small groups, or multi-clusters where velocity vectors remain relatively aligned.
\begin{figure}
    \centering
    \begin{subfigure}{0.24\textwidth}
        \centering
        \begin{overpic}
            [width=1.2\linewidth]{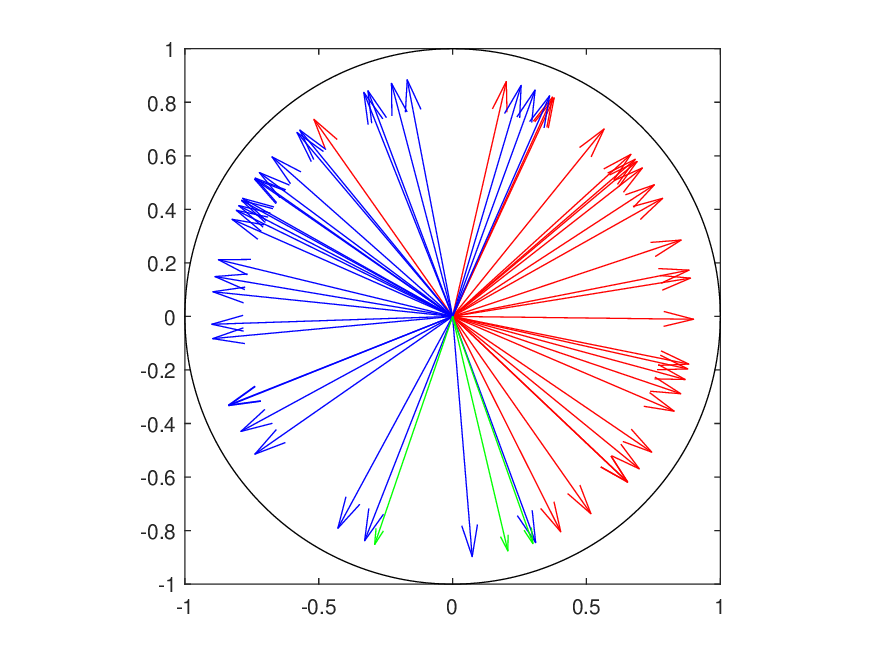}
            \put(37,73){\footnotesize (a) $t=0$}
        \end{overpic}
    \end{subfigure}
    \begin{subfigure}{0.24\textwidth}
        \centering
        \begin{overpic}
            [width=1.2\linewidth]{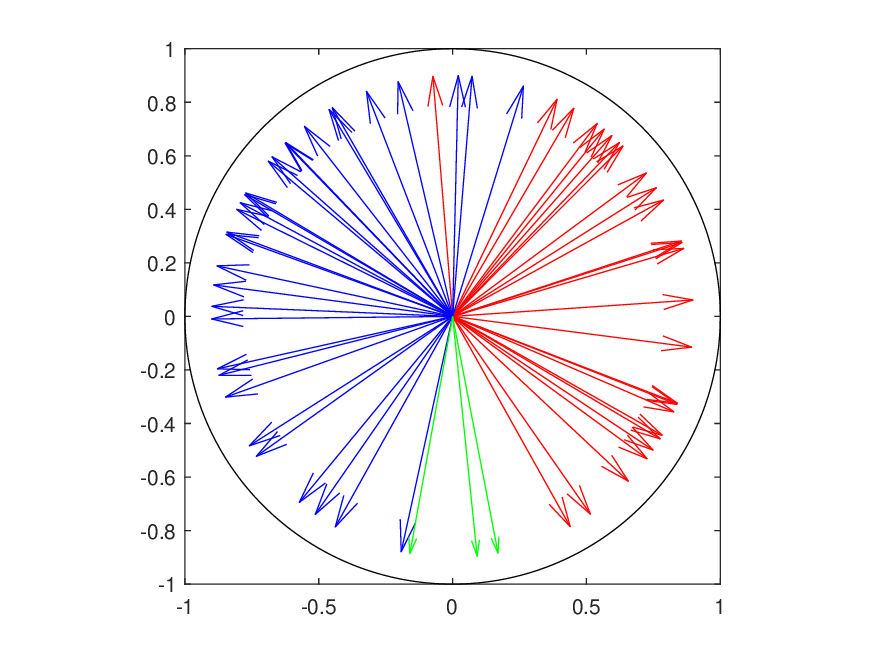}
            \put(37,73){\footnotesize (b) $t=5$}
        \end{overpic}
    \end{subfigure}
    \begin{subfigure}{0.24\textwidth}
        \centering
        \begin{overpic}
            [width=1.2\linewidth]{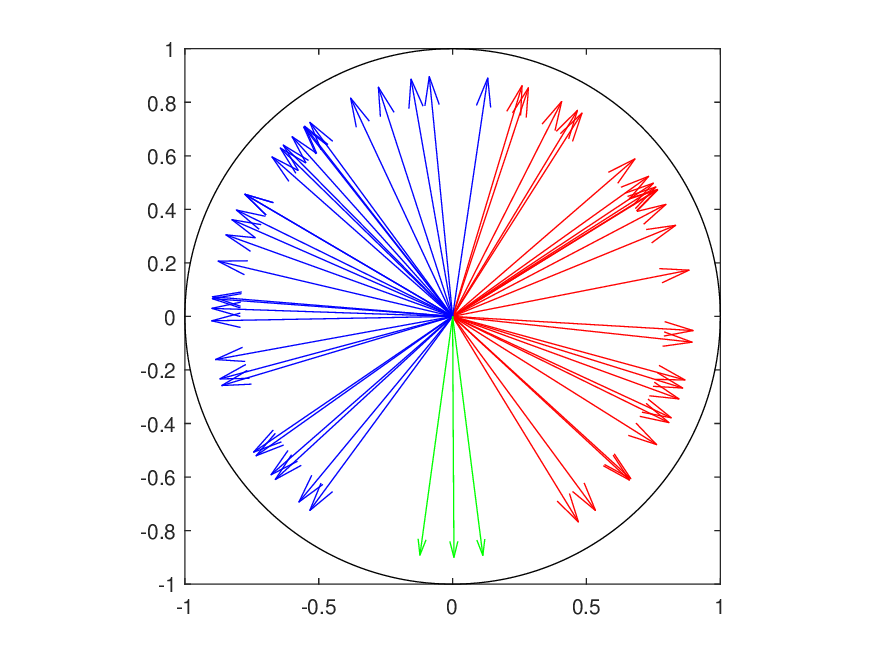}
            \put(36,73){\footnotesize (c) $t=10$}
        \end{overpic}
    \end{subfigure}
    \begin{subfigure}{0.24\textwidth}
        \centering
        \begin{overpic}
            [width=1.2\linewidth]{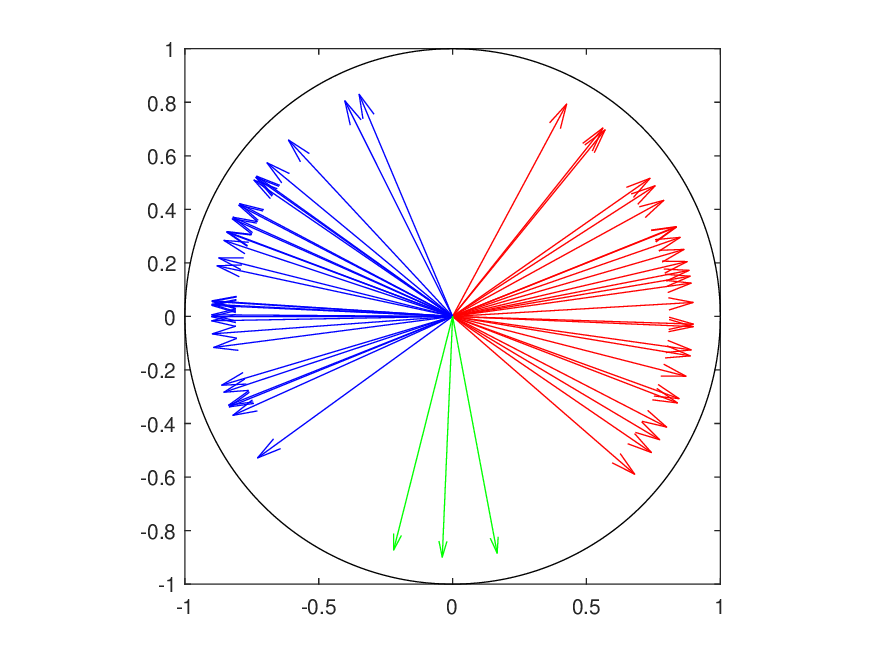}
            \put(36,73){\footnotesize (d) $t=30$}
        \end{overpic}
    \end{subfigure}
    \begin{subfigure}{0.24\textwidth}
        \centering
        \begin{overpic}
            [width=1.2\linewidth]{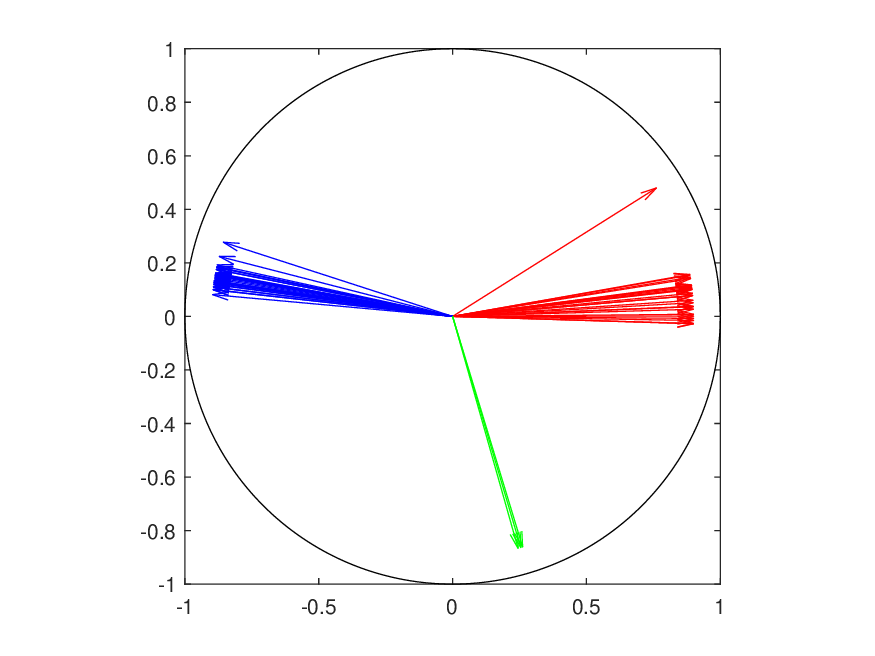}
            \put(36,73){\footnotesize (e) $t=50$}
        \end{overpic}
    \end{subfigure}
    \begin{subfigure}{0.24\textwidth}
        \centering
        \begin{overpic}
            [width=1.2\linewidth]{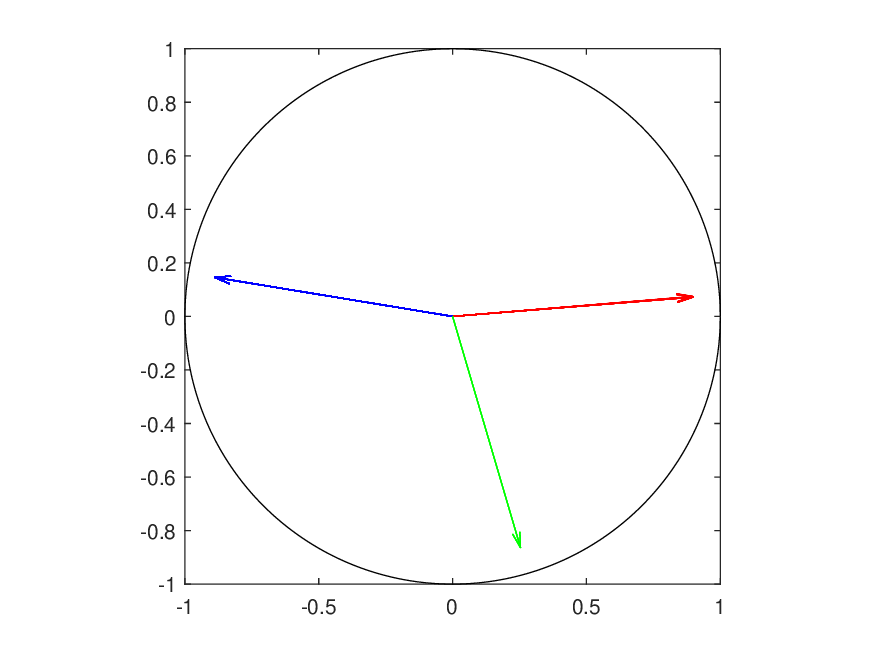}
            \put(36,73){\footnotesize (f) $t=80$}
        \end{overpic}
    \end{subfigure}
    \begin{subfigure}{0.24\textwidth}
        \centering
        \begin{overpic}
            [width=1.2\linewidth]{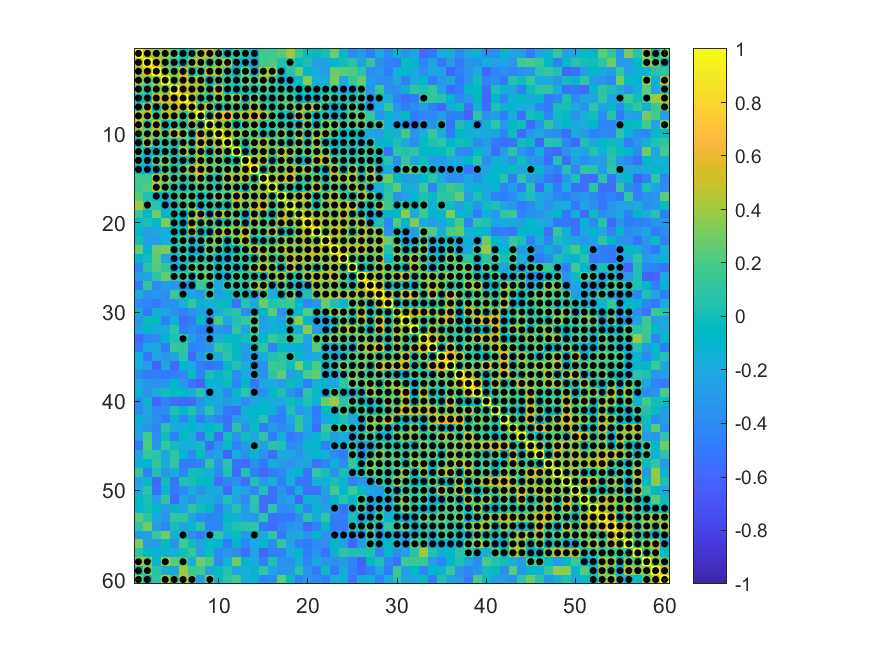}
            \put(33,73){\footnotesize (g) $t=30$}
        \end{overpic}
    \end{subfigure}
    \begin{subfigure}{0.24\textwidth}
        \centering
        \begin{overpic}
            [width=1.2\linewidth]{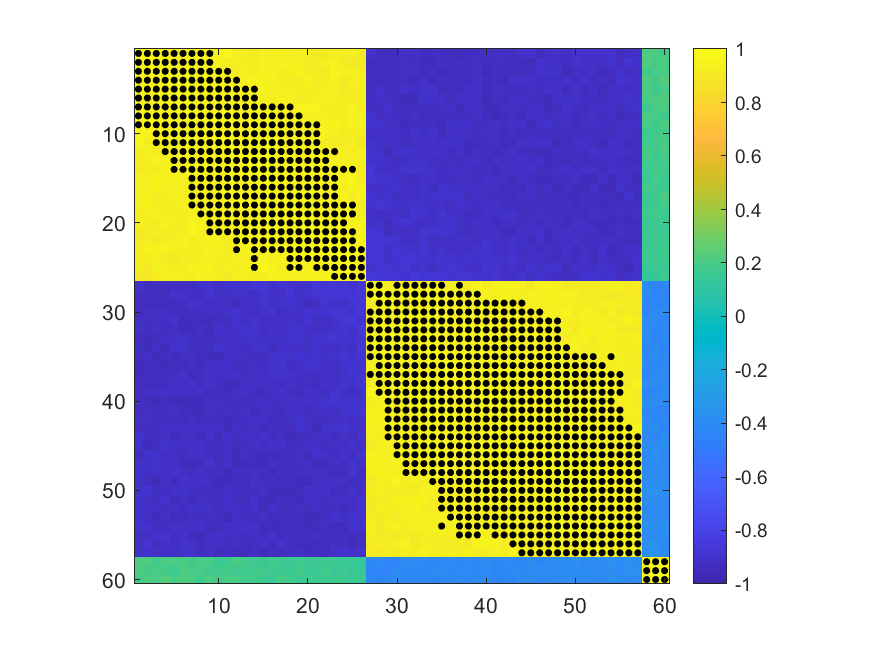}
            \put(32,73){\footnotesize (h) $t=300$}
            \put(90,47){\rotatebox{270}{\footnotesize $\kappa_{ij}$}}
        \end{overpic}
    \end{subfigure}
    \caption{Evolution over time for Case 6 in \Cref{tab_adapCS}. (a)-(f) show the angles of the normalized velocity vectors for all particles at six selected time-stamps, providing a snapshot of the system's configuration at each time. (g)-(h) display the values of $(\kappa_{ij})$ at $t=30$ and $t=300$, respectively. The values of $\kappa_{ij}$ are represented using MATLAB's \texttt{imagesc} function, where color indicates the magnitude of $\kappa_{ij}$.}
    \label{fig_snapshot_rng6}
\end{figure}
\begin{figure}
    \centering
    \begin{subfigure}{0.24\textwidth}
        \centering
        \begin{overpic}
            [width=1.2\linewidth]{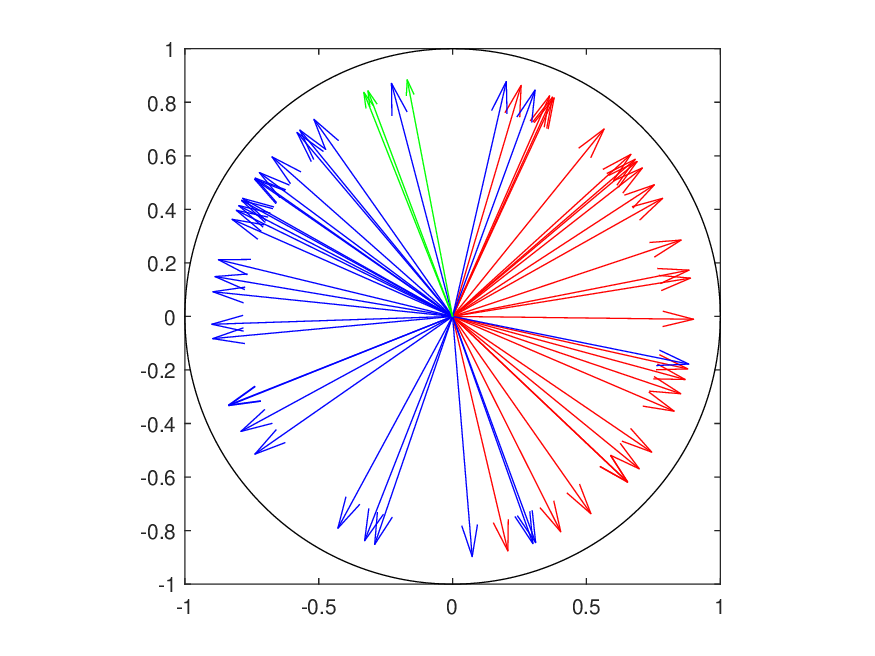}
            \put(37,73){\footnotesize (a) $t=0$}
        \end{overpic}
    \end{subfigure}
    \begin{subfigure}{0.24\textwidth}
        \centering
        \begin{overpic}
            [width=1.2\linewidth]{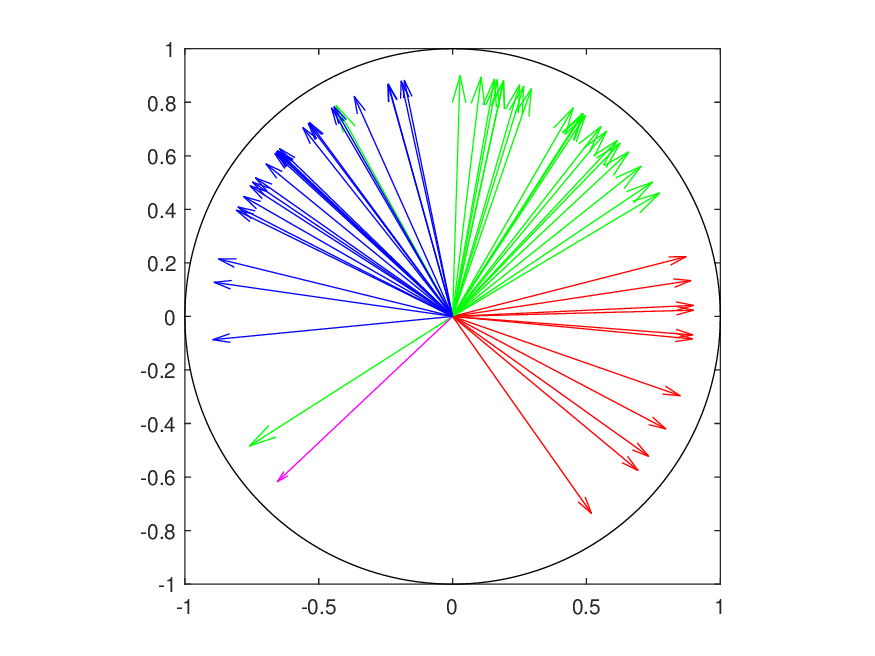}
            \put(37,73){\footnotesize (b) $t=2$}
        \end{overpic}
    \end{subfigure}
    \begin{subfigure}{0.24\textwidth}
        \centering
        \begin{overpic}
            [width=1.2\linewidth]{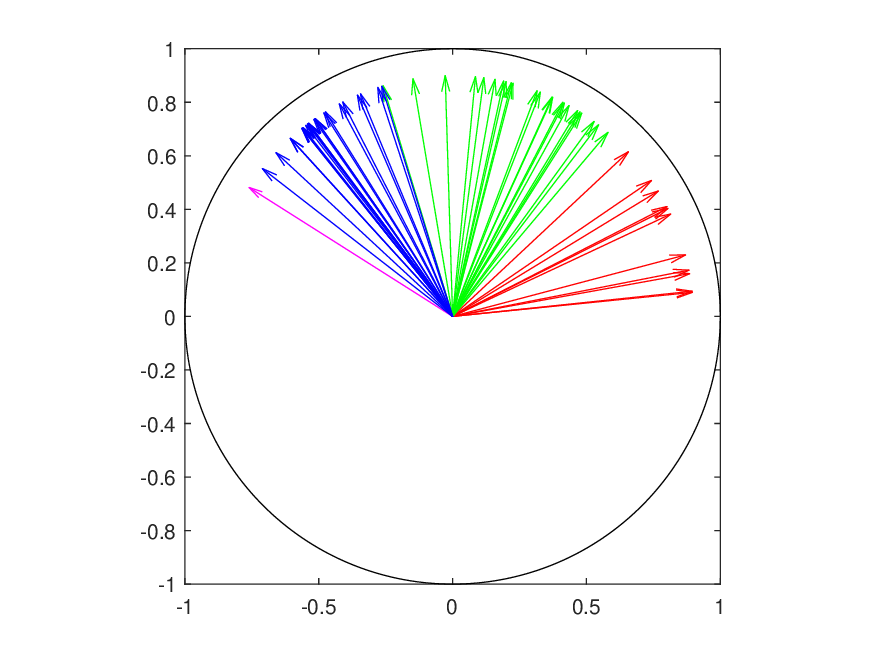}
            \put(37,73){\footnotesize (c) $t=5$}
        \end{overpic}
    \end{subfigure}
    \begin{subfigure}{0.24\textwidth}
        \centering
        \begin{overpic}
            [width=1.2\linewidth]{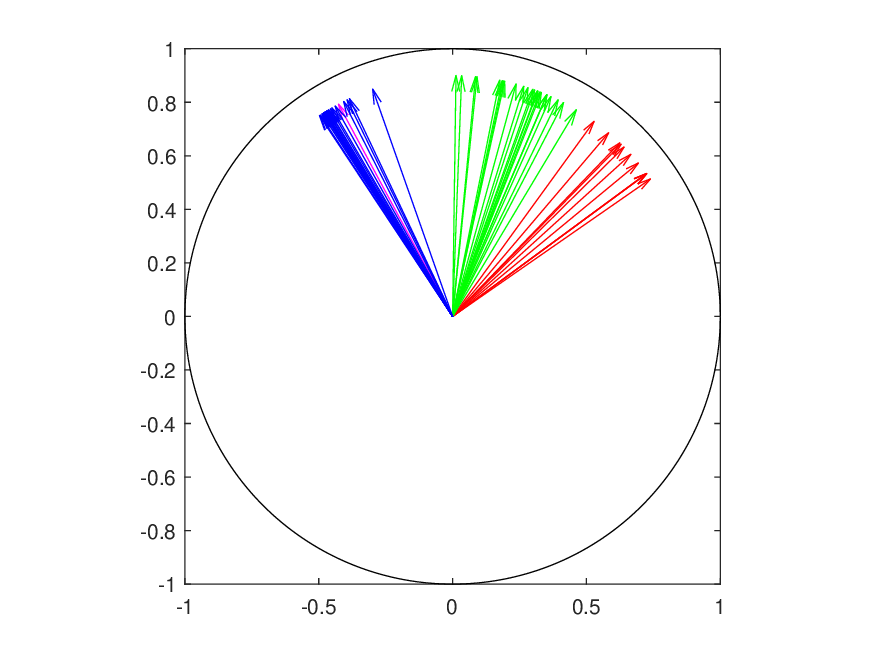}
            \put(36,73){\footnotesize (d) $t=10$}
        \end{overpic}
    \end{subfigure}
    \begin{subfigure}{0.24\textwidth}
        \centering
        \begin{overpic}
            [width=1.2\linewidth]{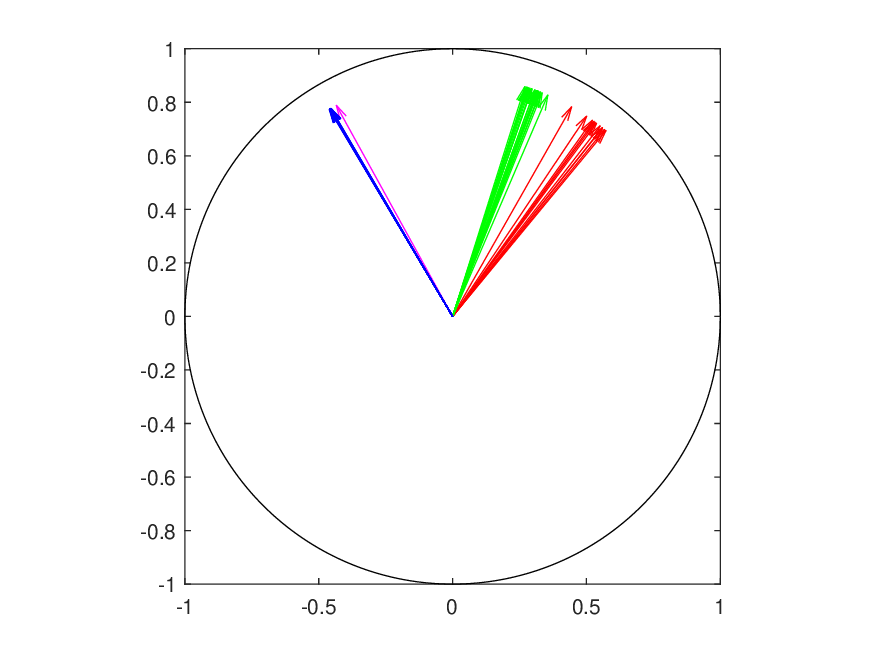}
            \put(36,73){\footnotesize (e) $t=20$}
        \end{overpic}
    \end{subfigure}
    \begin{subfigure}{0.24\textwidth}
        \centering
        \begin{overpic}
            [width=1.2\linewidth]{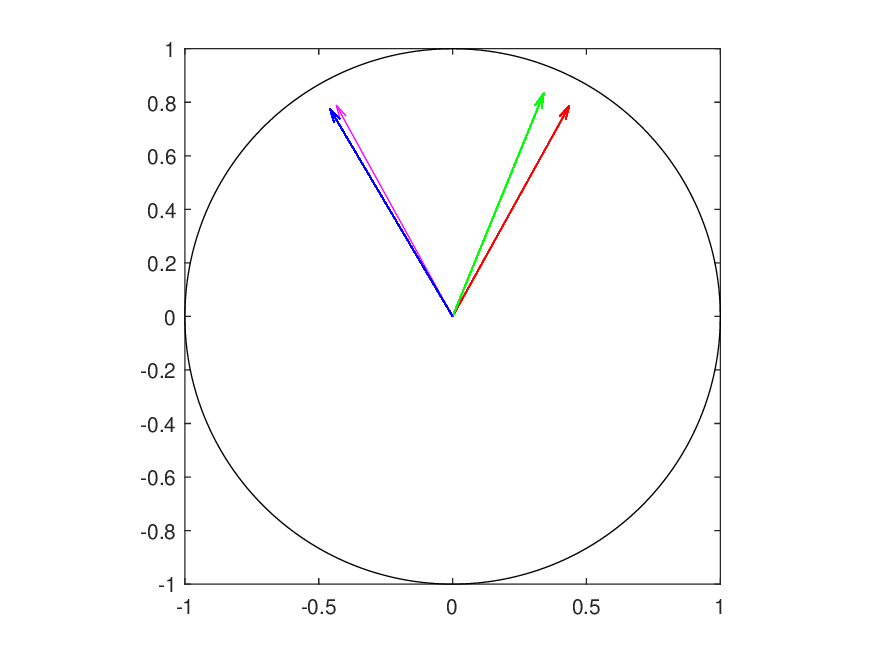}
            \put(36,73){\footnotesize (f) $t=80$}
        \end{overpic}
    \end{subfigure}
    \begin{subfigure}{0.24\textwidth}
        \centering
        \begin{overpic}
            [width=1.2\linewidth]{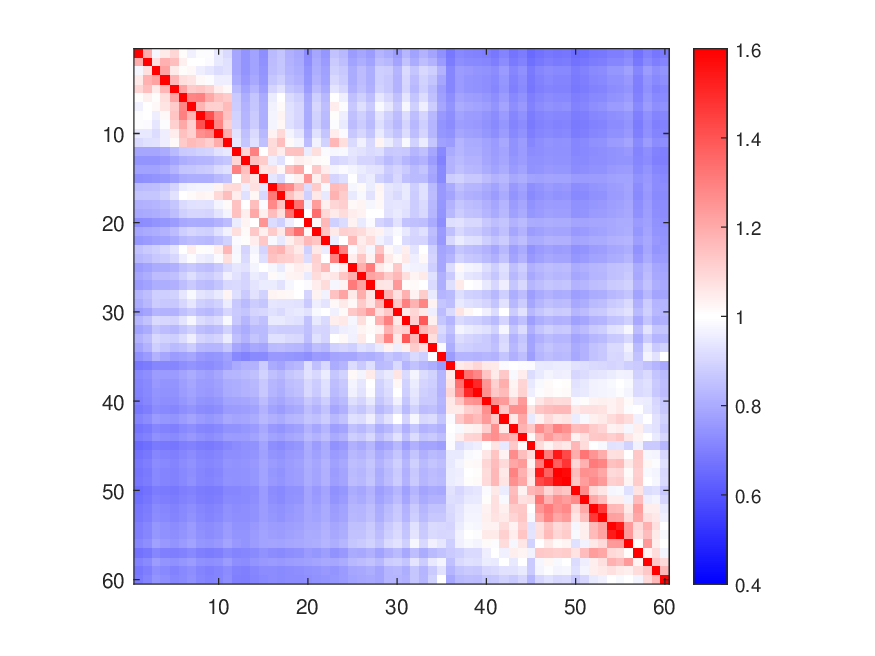}
            \put(33,73){\footnotesize (g) $t=10$}
        \end{overpic}
    \end{subfigure}
    \begin{subfigure}{0.24\textwidth}
        \centering
        \begin{overpic}
            [width=1.2\linewidth]{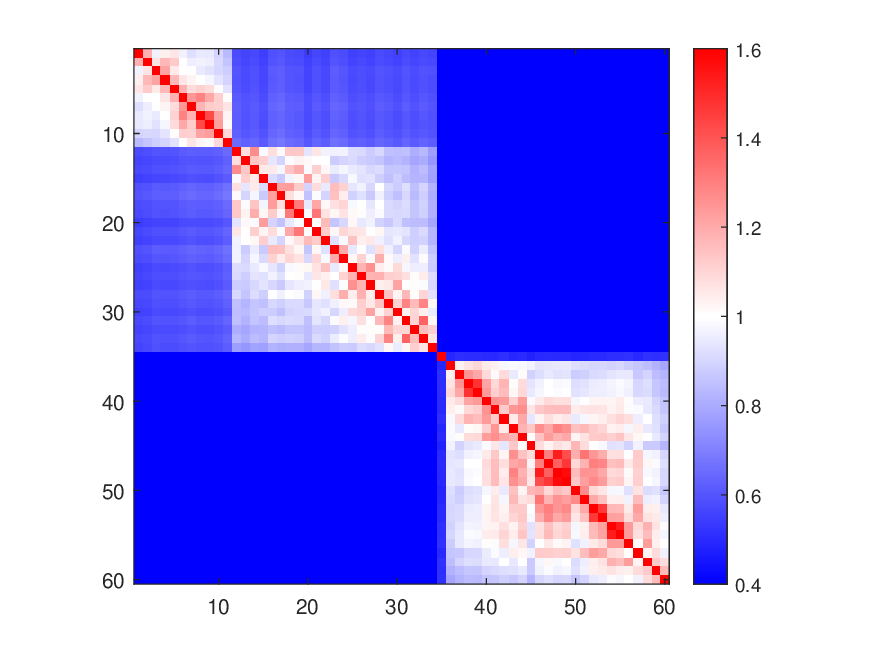}
            \put(32,73){\footnotesize (h) $t=300$}
            \put(90,61){\rotatebox{270}{\footnotesize $\log_d\|\bx_i-\bx_j\|$}}
        \end{overpic}
    \end{subfigure}
    \caption{Evolution over time for Case 6 in \Cref{tab_CS} when $d=0.009$. (a)-(f) show the angles of the normalized velocity vectors for all particles at six selected time-stamps, providing a snapshot of the system's configuration at each time. (g)-(h) display the spatial distance $\|\bx_i-\bx_j\|$, using the logarithm with base $d$, at $t=10$ and $t=300$, respectively. The values are represented using MATLAB's \texttt{imagesc} function, where red and blue indicates connectivity and disconnection, respectively.}
    \label{fig_snapshot_rng6_orig}
\end{figure}

To analyze the clustering process and opinion divergence, and to highlight the differences between the adaptive and classical CS dynamics, we present snapshots of the velocity vectors' polar angles in \Cref{fig_snapshot_rng6} and \Cref{fig_snapshot_rng6_orig} for Case 6. The selected time stamps capture key stages of the system's temporal evolution. Colors correspond to asymptotic group formation, while black dots indicate connections determined by the connection kernel $\psi$, illustrating current cluster formations.

First, consider the adaptive CS model in \Cref{fig_snapshot_rng6}. Initially, velocity polar angles exhibit somewhat vague groupings (see (a)). By $t=80$, multi-clusters are fully formed and stable (see (f)). To further support this, we plot the network variables $\kappa_{ij}$, which tends toward the similarity of opinions between $i$-th and $j$-th particles, at $t=30,300$ in (g)-(h). The similarity with value close to 1 (white) indicates agreement and near $-1$ (purple) indicates opposition. At $t=30$, the directional alignment of opinions became distinctly separated across groups, and synchronization within each group also began to emerge significantly (see (d)). In contrast, (g) reveals that the connectivity determined by determined by spatial distance through $\psi$, differs from the asymptotic cluster formation. This clearly suggests that the grouping of opinions—equivalent to the asymptotic group formation—occurs prior to the connection structure determined by the spatial kernel.

Next, examine the classical CS model in \Cref{fig_snapshot_rng6_orig}. At early times $(t\le10)$, all opinions converge within a half-circle despite initially distributed across the full circle. This inevitably implies that more than half of the opinions are disregarded in favor of a single dominant trend. Furthermore, multi-clusters emerge due to spatial disconnection. In (g)-(h), we plot the spatial distance $\|\bx_i-\bx_j\|$ on a logarithmic scale with base $d=0.009<1$, making $\log_d$ a decreasing function. Values equal to 1 (the connectivity threshold) appear white, values greater than 1 (connected) are red, and less than 1 (disconnected) are blue with lightness. By $t=10$, spatial cluster formation is nearly complete (cf. $t=300$), while opinions remain relatively dispersed compared to full synchronization. This indicates that multi-cluster formation in the classical CS model is driven purely by spatial separation in the process of reaching global opinion consensus.

These results highlight the adaptive CS model's ability to capture multi-cluster formation driven by opinion divergence, in contrast to the classical CS model where spatial disconnection dominates. This suggests that the adaptive CS model provides a more accurate framework for explaining opinion polarization in dynamic systems.

\subsection{Dynamics in the singular adaptive Cucker-Smale dynamics}
In this subsection, we analyze the effect of the radius parameter $d$
on group formation in the singular adaptive CS model \eqref{singular_adapCS}. All simulations are initialized with identical conditions as follows: with $N=50$,
\begin{align*}
    \theta_i=\frac{(i-1)\pi}{N-1},\quad\bx_i^0=(\cos\theta_i,\sin\theta_i),\quad\bv_i^0=-\frac{1}{3}(\cos\theta_i,\sin\theta_i),\quad\forall~i\in[N].
\end{align*}
We conducted simulations of the singular adaptive CS dynamics for various radius values within the range $d\in[0.5,0.65]$.
\begin{figure}
    \centering
    \begin{subfigure}{0.45\textwidth}
        \centering
        \begin{overpic}[width=1\linewidth]{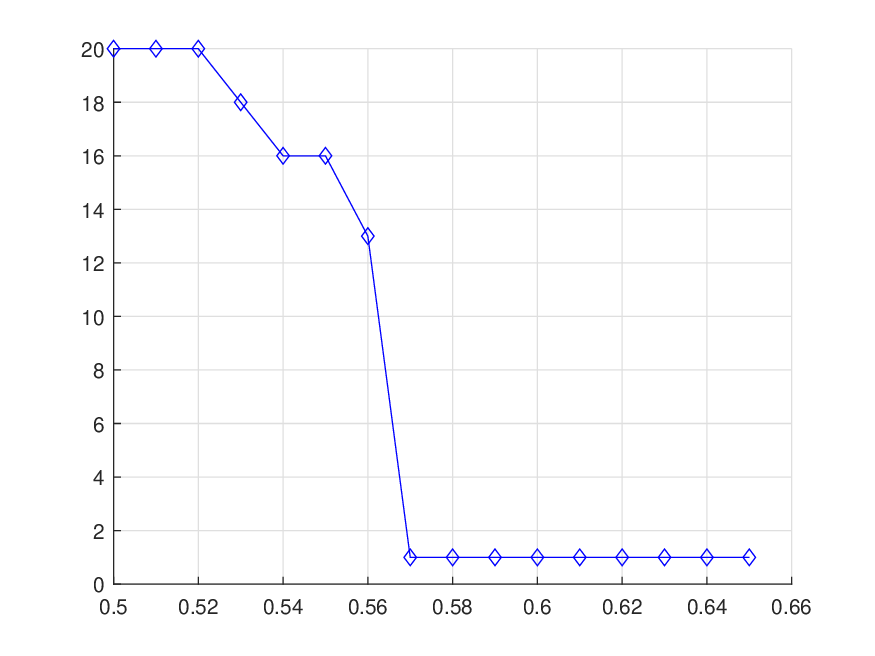}
            \put(3,17){\rotatebox{90}{\small Number of clusters}}
            \put(43, -1){\small radius $d$}
        \end{overpic}
        \caption{Number of clustering at $t=100$}
        \label{subfig_nbr_clusters}
    \end{subfigure}
    \begin{subfigure}{0.45\textwidth}
        \centering
        \begin{overpic}[width=1\linewidth]{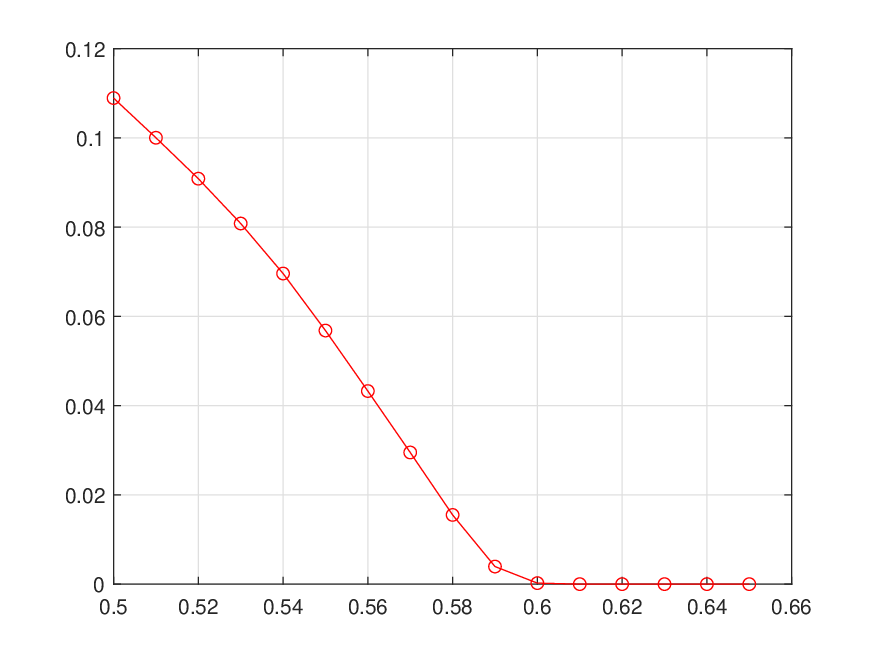}
            \put(1,19){\rotatebox{90}{\small mean of $\|\bv_i-\bv_c\|$}}
            \put(43, -1){\small radius $d$}
        \end{overpic}
        \caption{mean velocity deviation at $t=30$}
        \label{subfig_vel_dv}
    \end{subfigure}
    \caption{(Singular adaptive CS model) Result for different values of the radius parameter $d$. (a) The number of clusters obtained asymptotically at \(t = 100\) for various values of $d$. (b) The mean velocity deviation at $t = 30$ for different values of $d$. These illustrate how varying $d$ affects the system's clustering.}
    \label{sim_thm23}
\end{figure}
In \Cref{subfig_nbr_clusters}, the number of clusters at $t=100$ is depicted. This time is empirically validated as sufficient to observe asymptotic clustering behavior in this context. Here, clusters are defined as connected components in the temporal graph \eqref{CS_graph}. The results suggest the existence of a threshold radius $d$ for a given initial configuration that ensures the formation of a mono-cluster. As expected, the plot exhibits monotonic decrease, converging to 1 (mono-cluster).

To examine the influence of intermediate velocity configurations on asymptotic cluster structures, we plot the velocity deviation, defined by
\begin{align*}
    \frac{1}{N}\sum_{i=1}^N\|\bv_i-\bv_c\|\quad\left(\mbox{where}\quad\bv_c=\frac{1}{N}\sum_{i=1}^N\bv_i\right)
\end{align*}
at $t=30$ for the same set of radius values in \Cref{subfig_vel_dv}. When $d\ge0.6$, there is {\it almost} no velocity deviation among particles at $t=30$, indicating {\it near}-complete synchronization {\it already} at this time. However, \Cref{subfig_nbr_clusters} reveals that a mono-cluster emerges when $d\ge0.57$. In cases where $0.57\le d<0.6$, the positive velocity deviation at $t=30$ suggests differences of opinion still persist within subgroups, leading to {\it direct} disconnections, although graph connectivity is maintained through paths of length greater than $1$.
\begin{figure}
    \centering
    \begin{subfigure}{0.45\textwidth}
        \centering
        \begin{overpic}
            [width=0.8\linewidth]{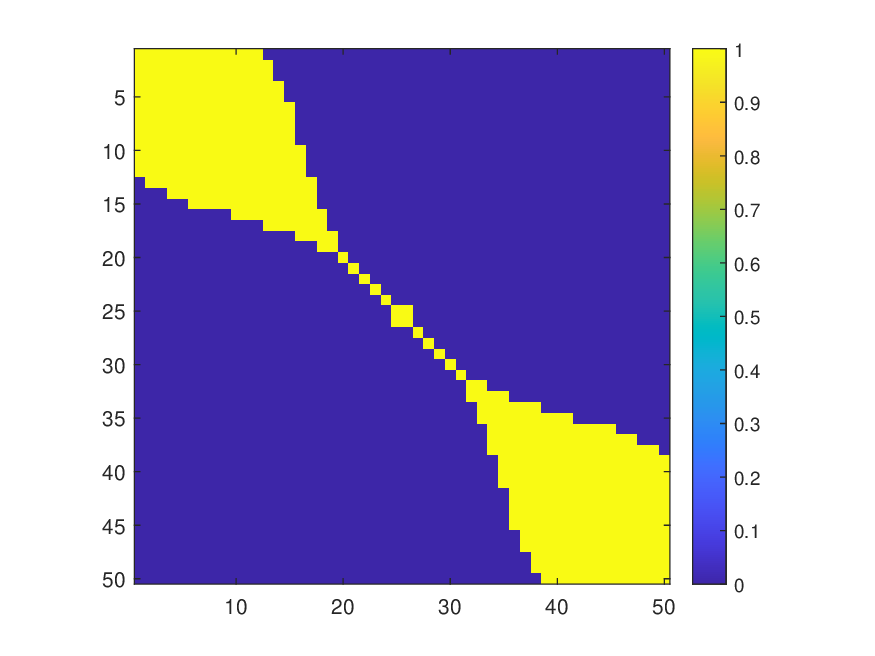}
            \put(30,73){\small (a) $d=0.56$}
            \put(90,47){\rotatebox{270}{\small $\psi_{ij}a_{ij}$}}
        \end{overpic}
    \end{subfigure}
    \begin{subfigure}{0.45\textwidth}
        \centering
        \begin{overpic}
            [width=0.8\linewidth]{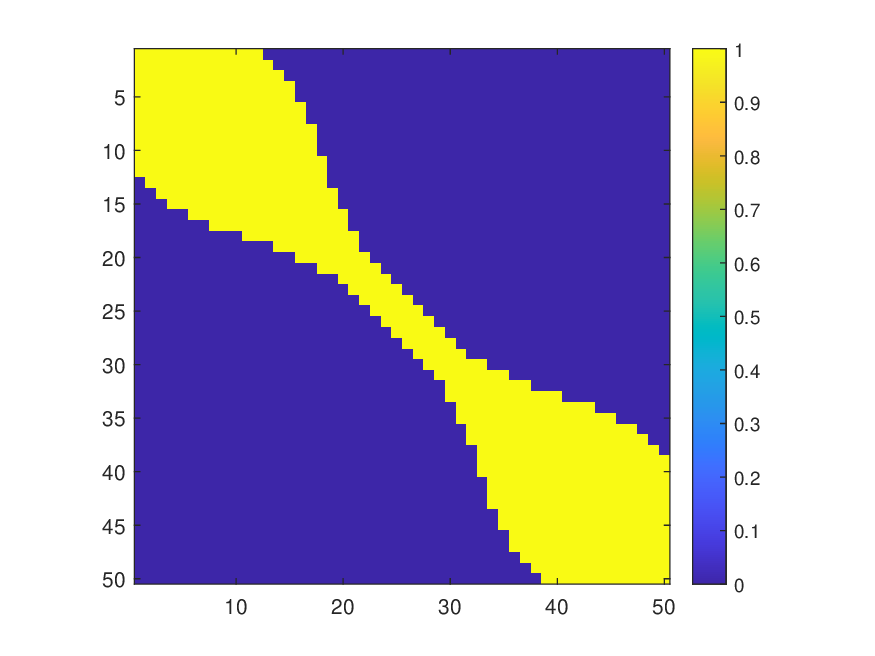}
            \put(30,73){\small (b) $d=0.58$}
        \end{overpic}
    \end{subfigure}
    \begin{subfigure}{0.45\textwidth}
        \centering
        \begin{overpic}
            [width=0.8\linewidth]{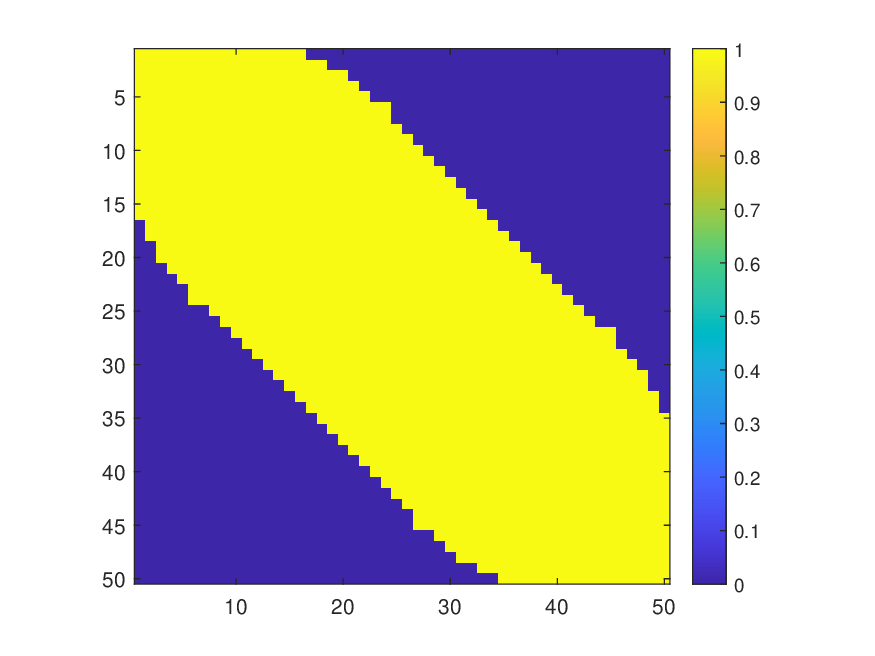}
            \put(30,73){\small (c) $d=0.60$}
        \end{overpic}
    \end{subfigure}
    \begin{subfigure}{0.45\textwidth}
        \centering
        \begin{overpic}
            [width=0.8\linewidth]{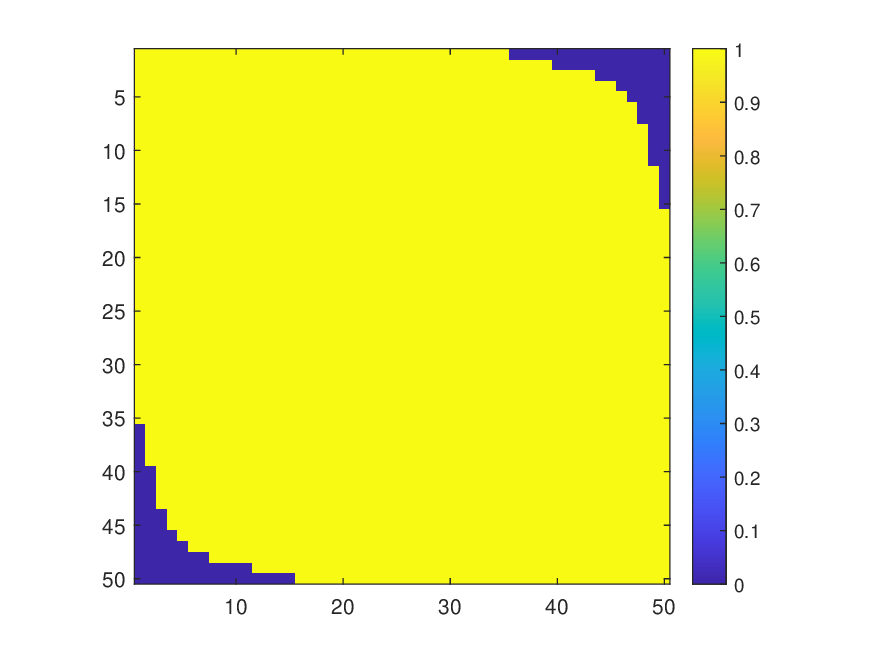}
            \put(30,73){\small (d) $d=0.65$}
        \end{overpic}
    \end{subfigure}
    \caption{Asymptotic graph structures for selected values of $d=0.56,0.58,0.6,0.65$ at $t=100$ from the results in \Cref{subfig_nbr_clusters}. The images show the connectivity and weight $\psi_{ij}a_{ij}$ of the graph for each selected $d$.}
    \label{sim_thm23_four}
\end{figure}

To illustrate this behavior, we consider four cases with radius values $d=0.56,0.58,0.6,0.65$. The results, shown in \Cref{sim_thm23_four}, depict the asymptotic ($t=100$) opinion similarity multiplied with connection kernel, i.e., $\psi_{ij}a_{ij}$. The clear yellow color in the plots illustrates complete synchronization within subgroups. For $d=0.56$, the system forms three distinct clusters along with ten isolated particles due to several disconnections. As $d$ increases, particles progressively form a single cohesive group, with connectivity strengthening through additional pairwise links. These results highlight that even within mono-clusters formed for different radius values $d$, the graph connectivity structure can vary significantly.

\section{Conclusion}\label{sec_concl}
In this paper, we proposed an adaptive rule for the Cucker-Smale (CS) model, specifically designed to describe the emergence of diverse opinions in opinion dynamics. From a mathematical perspective, we have focused on the singular limit case for fast adaptation. Due to Fenichel's Theorem, our results for this singular limit transfer to the full model. The singular limit model was formulated as Laplacian dynamics on temporal graphs, and we identified three classes of temporal graphs that guarantee flocking behavior. For each class, illustrative examples were provided to clarify their structure and applicability.

Additionally, we presented numerical experiments that demonstrate convergence in Laplacian dynamics for these graph classes. Through comparisons with the classical CS model, we highlighted the fully adaptive CS model's ability to capture complex opinion dynamics, including the polarization and formation of distinct clusters.

Our findings open several promising directions for future research, including the exploration of multi-clustering phenomena, the impact of graph topologies on dynamics, and potential applications to real-world systems.\medskip

\textbf{Acknowledgments:} J. Yoon would like to thank the Alexander von Humboldt Stiftung for support via a postdoctoral research fellowship.

\bibliographystyle{abbrv}

\end{document}